\documentclass[journal ]{new-aiaa}

\usepackage{graphicx}
\usepackage{amsmath}

\usepackage{enumitem}

\usepackage{caption}
\usepackage{float}

\usepackage{multirow}

\usepackage{bm}

\title{Analytical Strategies and Winning Conditions for Elliptic-Orbit Target-Attacker-Defender Game}

\author{Shuyue Fu \footnote{PhD Candidate, Shen Yuan Honors College, School of Astronautics, fushuyue@buaa.edu.cn.}, Shengping Gong \footnote{Professor, School of Astronautics, gongsp@buaa.edu.cn, Senior Member AIAA (Corresponding Author).}, Di Wu \footnote{Associate Professor, School of Astronautics, 1522620129@qq.com, Member AIAA.}, and Peng Shi\footnote{Professor, School of Astronautics, shipeng@buaa.edu.cn.}}
\affil{Beihang University, Beijing, 100191, People's Republic of China}

\begin{document}

\maketitle

\begin{abstract}
This paper proposes an analytical framework for the orbital Target-Attacker-Defender game with a non-maneuvering target along elliptic orbits. Focusing on the linear quadratic game, we derive an analytical solution to the matrix Riccati equation, which yields analytical Nash-equilibrium strategies for the game. Based on the analytical strategies, we derive the analytical form of the necessary and sufficient winning conditions for the attacker. The simulation results show good consistency between the analytical and numerical methods, exhibiting 0.004$\%$ relative error in the cost function. The analytical method achieves over 99.9$\%$ reduction in CPU time compared to the conventional numerical method, strengthening the advantage of developing the analytical strategies. Furthermore, we verify the proposed winning conditions and investigate the effects of eccentricity on the game outcomes. Our analysis reveals that for games with hovering initial states, the initial position of the defender should be constrained inside a mathematically definable set to ensure that the attacker wins the game. This constrained set further permits geometric interpretation through our proposed method. This work establishes the analytical framework for orbital Target-Attacker-Defender games, providing fundamental insights into the solution analysis of the game.
\end{abstract}

\section{Introduction}
\lettrine{D}{ifferential} game theory, first introduced by Isaacs \cite{isaacs1999differential}, has attracted significant attention from scholars. The differential game originally described the scenario involving two players, namely, the pursuer and evader \cite{pang2024solving}. In this scenario, the pursuer aims to intercept while the evader attempts to avoid capture \cite{jagat2017nonlinear,fu2025analytical}. The differential game theory has been widely applied in guidance of missiles \cite{liang2020optimal,Battistini6978887}, unmanned aerial vehicles \cite{Sun2017Multiple,wei2018Optimal}, spacecraft \cite{pang2024solving,li2019dimension}, and other autonomous vehicles \cite{prokopov2013linear,Li5751240}. Among these applications, the differential game applied in spacecraft guidance, i.e., the orbital differential game, presents particular challenges due to the strong dynamical constraints \cite{pang2024solving,li2024nash} and the high dimensionality of the problem \cite{li2019dimension}. To solve the orbital two-player game, several methods have been proposed. Tartaglia and Innocenti \cite{tartaglia2016game} investigated the orbital two-player game using the linear quadratic (LQ) method in the context of the Clohessy-Wiltshire (CW) equations along the circular reference orbits \cite{clohessy1960terminal}. Li et al. \cite{li2019dimension,li2020saddle} developed a dimension-reduction method and explored the Nash-equilibrium solutions for the orbital two-player games. Pang et al. \cite{pang2024solving} provided precise gradient to solve the two-player game along the elliptic reference orbits, in the context of the Tschauner-Hempel (TH) equations \cite{dang2017solutions}, and Fu et al. \cite{fu2025analytical} developed the analytical LQ strategies for the two-player game along arbitrary Keplerian orbits. These methods provide foundations for extending to multi-player scenarios \cite{chen2024luring,li2024nash,wu2025exterior}, such as the Target-Attacker-Defender (TAD) games \cite{li2024nash}.

In the TAD game, there are three players, namely, the attacker, the defender, and the target. In this scenario, the attacker attempts to capture the target, and the defender works to intercept the attacker. Recently, this scenario has gained increasing attention in space security applications \cite{qian2025swarm}. The design of control strategies in the TAD games proves more complex than in the two-player scenarios due to the interception role of the defender \cite{Zhouapp9153190}. For this consideration, Liu et al. \cite{liu2018Optimal} developed a hierarchical optimization method to solve the TAD game using the impulsive maneuver. When considering the game with continuous thrust, Bian and Zheng \cite{bian2023control} and Li \cite{li2024orbital} developed the LQ method for the TAD game. Furthermore, Li et al. \cite{10.1007/978-981-97-3324-8_28,li2024nash} established the models of the TAD game with a non-maneuvering target, including the LQ game and the game scenarios with thrust constraints. For the LQ method, the solution of the problem can be reduced to solving the matrix Riccati equation. The Riccati equation can be categorized into the differential Riccati equation (DRE) and the algebraic Riccati equation (ARE). The DRE applies to finite-horizon problems while the ARE is suitable for infinite-horizon cases \cite{Kumar8818296}. For the orbital differential game, given the practical relevance of finite-horizon problems \cite{liao2021research,li2024orbital}, we focus on the finite-horizon problem and the solution of the DRE. In the following text, the matrix Riccati equation refers to the DRE. In previous works on the LQ TAD game \cite{li2024orbital,10.1007/978-981-97-3324-8_28,li2024nash}, the matrix Riccati equation was solved by conventional numerical methods, such as numerical backward integration and discretization method. However, scenarios requiring fast response times necessitate analytical solutions to improve computational efficiency. Our previous work \cite{fu2025analytical} provided an analytical solution of the matrix Riccati equation for the orbital two-player game, yielding the analytical LQ strategies. In this paper, we extend the analytical LQ strategies of the TAD game with a non-maneuvering target based on the analytical solution of the TH equations \cite{dang2017solutions} and variable transformation \cite{li2012fuel}, illustrating the extensibility of our developed method. 

For the TAD game with the fixed game time, the analytical strategies enable the analytical solutions because initial costate variables are analytically accessible. This leads to a critical research question: Under what conditions does the attacker win the game (i.e., what are the winning conditions for the attacker)? Due to the three players involved, the analysis of winning conditions for the TAD game proves more complex \cite{li2024nash}. For such a scenario, Li et al. \cite{10.1007/978-981-97-3324-8_28,li2024nash} provided a pioneering exploration of this question. They pointed out that under the assumption that the attacker and defender have the hovering initial states, i.e., zero initial velocity of the controlled CW or TH equations, the initial position of the attacker and defender determines the game outcomes. However, their analysis relied on the numerical results and lacked geometric interpretation. Our work extends this analysis by deriving the analytical form of necessary and sufficient winning conditions with a geometric interpretation. Specifically, when the initial position of the attacker is given, the initial position of the defender should be constrained inside an analytically expressible set. We further analyze how orbital eccentricity affects the game outcomes. Based on the aforementioned discussion, two main contributions of this work can be summarized as follows:

\begin{enumerate}[label=(\arabic*)]
\item We derive the analytical strategies using the LQ method for the TAD game with a non-maneuvering target along the elliptic orbits. When the game time is fixed, these analytical strategies yield the analytical solution of the TAD game because the initial costate variables can be analytically obtained. This analytical solution achieves over 99.9$\%$ CPU time reduction compared to the numerical method, which is suitable for real-time applications requiring fast response times.
\item Based on the obtained analytical strategies, the analytical sufficient and necessary winning conditions for the attacker are further derived, establishing defender position constraints as mathematically definable sets. Furthermore, we provide the geometric interpretation for these conditions.
\end{enumerate}

The rest of this paper is organized as follows. Section \ref{sec2} presents the problem statement, including the controlled TH equations and the definition of the LQ TAD game. Section \ref{sec3} derives the analytical LQ strategies along the elliptic orbits. Section \ref{sec4} presents the analytical sufficient and necessary conditions ensuring that the attacker wins based on the obtained analytical strategies. The effectiveness of the proposed analytical strategies and winning conditions are verified in Section \ref{sec5}. Finally, conclusions are drawn in Section \ref{sec6}.

\section{Problem Statement}\label{sec2}
In this section, we provide the theoretical foundation of the orbital TAD game with a non-maneuvering target, including the controlled Tschauner-Hempel (TH) equations for the elliptic orbits and the model of the considered TAD game.

\subsection{Controlled Tschauner-Hempel Equations}\label{subsec2.1}

\begin{figure}[h]
\centering
\includegraphics[width=0.4\textwidth]{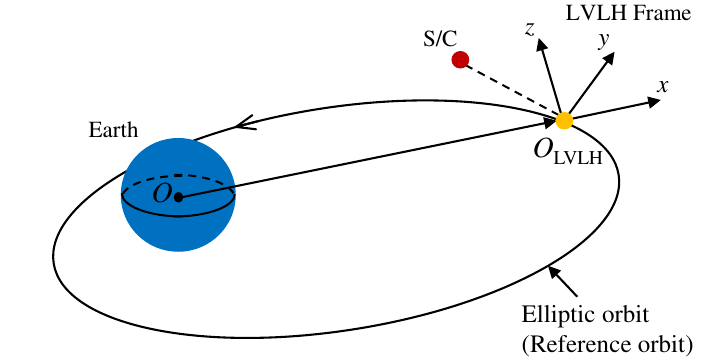}
\caption{The LVLH frame.}
\label{fig1}
\end{figure}

Since the game time is typically brief and the distance between the spacecraft is usually close \cite{li2024nash,fu2025analytical}, the problem can be simplified by using the linear relative dynamics \cite{stupik2012optimal} with negligible error \cite{gong9103281,gim2003state}. Furthermore, practical relevance makes it suitable to start games when proximity ensures mutual detection \cite{li2024nash,Zhouapp9153190}. Therefore, using this simplification, the orbital differential game can be solved by taking advantage of the analytical solution (i.e., state transition matrix (STM) \cite{yamanaka2002new,dang2017solutions,dang2017new}) of the linear relative dynamics \cite{dang2017solutions}. Since we focus on the TAD game along the elliptic orbits, the controlled TH equations \cite{li2012fuel} are introduced to describe the motion of the attacker, defender, and target. When expressing the controlled TH equations, the local-vertical/local-horizontal (LVLH) frame \cite{wen2016relative} is adopted, as shown in Fig. \ref{fig1}. 

The LVLH frame, rotating relative to the Earth-centered inertial (ECI) frame, originates at a reference orbit point. The \textit{x} axis extends radially outward, the \textit{z} axis aligns with orbital angular momentum ($\bm{h}$), and the \textit{y} axis is determined by the right-hand rule. Defining relative state as $\bm{X}=\left[x,\text{ }y,\text{ }z,\text{ }\dot{x},\text{ }\dot{y},\text{ }\dot{z}\right]^{\text{T}}$ (${\dot{ \left(  \cdot  \right)} }$ denotes the time derivative), the controlled TH equations are expressed as \cite{li2012fuel}:

\begin{equation}
\left\{ \begin{gathered}
  \ddot x = 2\dot f\dot y + \ddot fy + {\left( {\dot f} \right)^2}x + 2\mu {\left( {\frac{{\dot f}}{h}} \right)^{3/2}}x + {u_x} \hfill \\
  \ddot y =  - 2\dot f\dot x - \ddot fx + {\left( {\dot f} \right)^2}y - \mu {\left( {\frac{{\dot f}}{h}} \right)^{3/2}}y + {u_y} \hfill \\
  \ddot z =  - \mu {\left( {\frac{{\dot f}}{h}} \right)^{3/2}}z + {u_z} \hfill \\ 
\end{gathered}  \right.\label{eq100000}
\end{equation}
where $h$ denotes the value of $\bm{h}$, $\mu$ denotes the Earth gravitational constant, and $f$ denotes true anomaly of the reference orbit. The vector $\bm{u}=\left[u_x,\text{ }u_y,\text{ }u_z\right]^{\text{T}}$ is the vector of control inputs. To simplify the expression of the TH equations, the following transformation is adopted \cite{yamanaka2002new}:

\begin{equation}
\left[ {\begin{array}{*{20}{c}}
  {\tilde x} \\ 
  {\tilde y} \\ 
  {\tilde z} 
\end{array}} \right] = \rho \left[ {\begin{array}{*{20}{c}}
  x \\ 
  y \\ 
  z 
\end{array}} \right]\label{eq1}
\end{equation}
where $\rho=1+e\cos{f}$. The parameter $e$ denotes the eccentricity of the reference orbit. Using the relationship shown in Eq. \eqref{eq1}, we transform the controlled TH equations into the true anomaly domain (${\left(  \cdot  \right)^\prime }$ denotes the true anomaly derivative) \cite{li2012fuel}:

\begin{equation}
\left\{ \begin{gathered}
  \tilde x'' - \frac{3}{\rho }\tilde x - 2\tilde y' = \frac{\beta }{{{\rho ^3}}}{u_x} \hfill \\
  \tilde y'' + 2\tilde x' = \frac{\beta }{{{\rho ^3}}}{u_y} \hfill \\
  \tilde z'' + \tilde z = \frac{\beta }{{{\rho ^3}}}{u_z} \hfill \\ 
\end{gathered}  \right.\label{eq2}
\end{equation}
where
\begin{equation}
\beta  = \frac{1}{{{n^2}}} = \frac{1}{{{\mu  \mathord{\left/
 {\vphantom {\mu  {{p^3}}}} \right.
 \kern-\nulldelimiterspace} {{p^3}}}}}\label{eq2222222}
\end{equation}
where $p$ denotes the semilatus rectum of the reference orbit. For the elliptic orbits, $p=a\left(1-e^2\right)$, where $a$ denotes the semi-major axis of the reference orbit. Therefore, Eq. \eqref{eq2} can be expressed in the following matrix form:

\begin{equation}
\bm{\tilde X'} = \bm{A\tilde X} + \bm{Bu}\label{eq4}
\end{equation}
The matrices $\bm{A}$ and $\bm{B}$ are expressed as:

\begin{equation}
\bm{A} = \left[ {
  {\begin{array}{*{20}{c}}
  0&0&0&1&0&0 \\ 
  0&0&0&0&1&0 \\ 
  0&0&0&0&0&1 \\
  {\frac{3}{\rho }}&0&0&0&2&0 \\
  0&0&0&{ - 2}&0&0 \\
  0&0&{ - 1}&0&0&0
\end{array}}} \right]\label{eq5}
\end{equation}

\begin{equation}
\bm{B} = \frac{\beta }{{{\rho ^3}}}\left[ {\begin{array}{*{20}{c}}
  {{\bm{O}_{3 \times 3}}} \\ 
  {{\bm{I}_{3 \times 3}}} 
\end{array}} \right]\label{eq6}
\end{equation}

\subsection{TAD Game with A Non-Maneuvering Target}\label{subsec2.2}

\begin{figure}[h]
\centering
\includegraphics[width=0.4\textwidth]{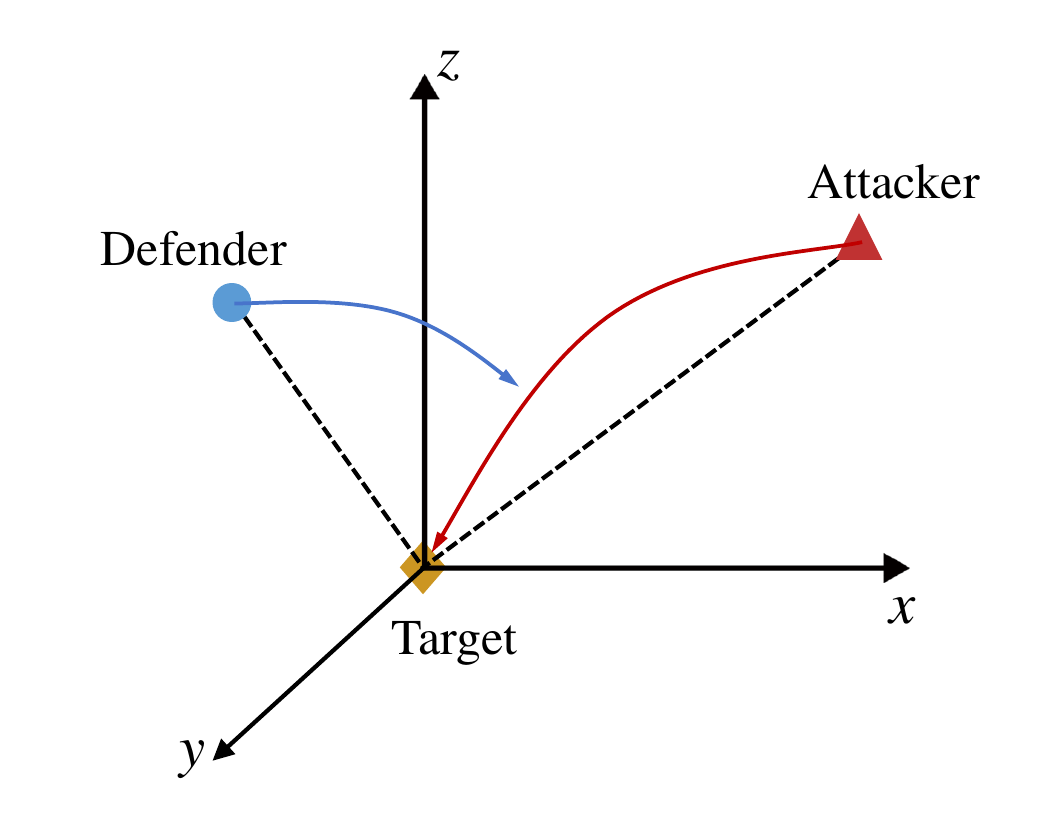}
\caption{The schematic of the TAD game in the LVLH frame.}
\label{fig2}
\end{figure}

In this subsection, the model of the TAD game with a non-maneuvering target is presented. This model was previously developed by Li et al. \cite{li2024nash} in terms of the two-dimensional CW relative dynamics (i.e., planar TAD game along the circular orbits). In this paper, we extend this model to the three-dimensional game scenario along elliptic orbits. The schematic of the TAD game is shown in Fig. \ref{fig2}. 

There are three players in the game, i.e., an attacker, a defender, and a target. The target has no maneuvering capability, the attacker attempts to capture the target, while the defender works to intercept the attacker. First, the following two assumptions of the considered TAD game are presented \cite{li2024nash}:

\begin{enumerate}[label=(\arabic*)]
\item The complete-information assumption is adopted, i.e., the state information about the attacker, the defender, and the target can be accurately obtained during the game.
\item The TAD game with a fixed game time is considered, i.e., a fixed interval of true anomaly along elliptic orbits.
\end{enumerate}

Then, we present the dynamical equations of the considered TAD game. The motion of the attacker, defender, and target follows the TH equations presented in Subsection \ref{subsec2.1} (note that the target has no maneuvering capability):
\begin{equation}
\left\{ \begin{gathered}
  {{\bm{\tilde X'}}_a} = \bm{A}{{\bm{\tilde X}}_a} + \bm{B}{\bm{u}_a} \hfill \\
  {{\bm{\tilde X'}}_d} = \bm{A}{{\bm{\tilde X}}_d} + \bm{B}{\bm{u}_d} \hfill \\ 
  {{\bm{\tilde X'}}_t} = \bm{A}{{\bm{\tilde X}}_t} \hfill \\
\end{gathered}  \right.\label{eq7}
\end{equation}
The subscript “\textit{a}” denotes the quantities corresponding to the attacker, the subscript “\textit{d}” denotes the quantities corresponding to the defender, and the subscript “\textit{t}” denotes the quantities corresponding to the target. Similar to Ref. \cite{li2024nash}, the target is assumed to be located at the origin of the LVLH frame, i.e., ${\bm{\tilde X}_t}=\bm{0}$. Therefore, defining $\bm{\tilde x}_{at} = {\bm{\tilde X}_a} - {\bm{\tilde X}_t}={\bm{\tilde X}_a}$ and $\bm{\tilde x}_{da} = {\bm{\tilde X}_d} - {\bm{\tilde X}_a}$, the dynamical equations of the TAD game can be rewritten as:

\begin{equation}
\left\{ \begin{gathered}
  {\bm{\tilde X'}_a} = \bm{A\tilde X}_{a} + \bm{B}{\bm{u}_a} \hfill \\
  \bm{\tilde x'}_{da} = \bm{A\tilde x}_{da} + \bm{B}{\bm{u}_d} - \bm{B}{\bm{u}_a} \hfill \\ 
\end{gathered}  \right.\label{eq8}
\end{equation}
Since the linear quadratic (LQ) method can generate a feedback and optimal control strategy, we employ the LQ method to generate the Nash-equilibrium strategies \cite{li2024nash,fu2025analytical}. The cost function of the LQ TAD game can be expressed as follows \cite{li2024nash}: 

\begin{equation}
J = \frac{1}{2}{\bm{\tilde X}_a}^{\text{T}}\left( {{f_f}} \right){\bm{S}_a}{\bm{\tilde X}_a}\left( {{f_f}} \right) - \frac{1}{2}{\bm{\tilde x}_{da}}^{\text{T}}\left( {{f_f}} \right){\bm{S}_d}_a{\bm{\tilde x}_{da}}\left( {{f_f}} \right) + \frac{1}{2}\int_{{f_0}}^{{f_f}} {\left( {{\bm{u}_a}^{\text{T}}{\bm{R}_a}{\bm{u}_a} - {\bm{u}_d}^{\text{T}}{\bm{R}_d}{\bm{u}_d}} \right)} {\text{d}}f\label{eq9}
\end{equation}
where $\bm{S}_a={\text{diag}}\left( s_{ar}{{\bm{I}_{3 \times 3}},{\text{ }}s_{av}{\bm{I}_{3 \times 3}}} \right)$, $\bm{S}_{da}={\text{diag}}\left( s_{dar}{{\bm{I}_{3 \times 3}},{\text{ }}s_{dav}{\bm{I}_{3 \times 3}}} \right)$, $\bm{R}_a={r_a}{\bm{I}_{3 \times 3}}$, and $\bm{R}_d={r_d}{\bm{I}_{3 \times 3}}$ denote the weighting matrices. Among these matrices, $\bm{S}_a$ and $\bm{S}_{da}$ are semi-positive definite matrices, while $\bm{R}_a$ and $\bm{R}_d$ are positive definite matrices. The $f_0$ and $f_f$ indicate the initial and terminal true anomaly of the reference orbit during the TAD game. The Hamiltonian function is defined as:

\begin{equation}
\mathcal{H} = \frac{1}{2}\left( {{\bm{u}_a}^{\text{T}}{\bm{R}_a}{\bm{u}_a} - {\bm{u}_d}^{\text{T}}{\bm{R}_d}{\bm{u}_d}} \right) + {\bm{\lambda }^{\text{T}}}\left( {\bm{A\tilde X}_a + \bm{B}{\bm{u}_a}} \right) + {{\bm{\nu }^{\text{T}}}\left( {\bm{A\tilde x}_{da} + \bm{B}{\bm{u}_d}-\bm{B}{\bm{u}_a}} \right)}\label{eq10}
\end{equation}
where $\bm{\lambda }$ and $\bm{\nu }$ denotes the costate variables. The transversality functional is defined as:
\begin{equation}
\Phi  = \frac{1}{2}{\bm{\tilde X}_a}^{\text{T}}\left( {{f_f}} \right){\bm{S}_a}{\bm{\tilde X}_a}\left( {{f_f}} \right) - \frac{1}{2}{\bm{\tilde x}_{da}}^{\text{T}}\left( {{f_f}} \right){\bm{S}_d}_a{\bm{\tilde x}_{da}}\left( {{f_f}} \right)\label{eq18}
\end{equation}
Let $\bm{\tilde y} = {\left[ {{{\bm{\tilde X}}^{\text{T}}_a},{\text{ }}{{\bm{\tilde x}}^{\text{T}}_{da}}} \right]^{\text{T}}}$ and $\bm{\Lambda } = {\left[ {\bm{\lambda }^{\text{T}},{\text{ }}\bm{\nu }}^{\text{T}} \right]^{\text{T}}}$, the transversality condition is expressed as:
\begin{equation}
\bm{\Lambda }\left( {{f_f}} \right) = \left.\frac{{\partial \Phi }}{{\partial \bm{\tilde y}}}\right|_{\bm{\tilde y}=\bm{\tilde y}\left(f_f\right)} = \left[ {\begin{array}{*{20}{c}}
  {{\bm{S}_a}}&{{\bm{O}_{6 \times 6}}} \\ 
  {{\bm{O}_{6 \times 6}}}&{-{\bm{S}_{da}}} 
\end{array}} \right]\bm{\tilde y}\left( {{f_f}} \right)
\label{eq19}
\end{equation}
The Nash-equilibrium strategies (i.e., ${\bm{u}_a}^*$ and ${\bm{u}_d}^*$) of the considered TAD game satisfy \cite{li2024nash}:

\begin{equation}
\mathop {\min }\limits_{{\bm{u}_a}^*} \mathop {\max }\limits_{{\bm{u}_d}^*} J\label{eq1222}
\end{equation}
According to Ref. \cite{li2024nash}, the Nash-equilibrium strategies of the TAD game are presented as follows:

\begin{equation}
\left\{ \begin{gathered}
  \left.\frac{{\partial \mathcal{H}}}{{\partial {\bm{u}_a}}}\right|_{{{\bm{u}_a} = {\bm{u}_a}^*}} = 0 \Rightarrow {\bm{u}_a}^* =  - {\bm{R}_a}^{ - 1}{\bm{B}^{\text{T}}}\left( {\bm{\lambda }\left( f \right) - \bm{\nu }\left( f \right)} \right) \hfill \\
  \left.\frac{{\partial \mathcal{H}}}{{\partial {\bm{u}_d}}}\right|_{{{\bm{u}_d} = {\bm{u}_d}^*}} = 0 \Rightarrow {\bm{u}_d}^* = {\bm{R}_d}^{ - 1}{\bm{B}^{\text{T}}}\bm{\nu }\left( f \right) \hfill \\ 
\end{gathered}  \right.
\label{eq12}
\end{equation}

Substituting these strategies into Eq. \eqref{eq8}, we obtain the following equations:

\begin{equation}
\left\{ \begin{gathered}
  {\bm{\tilde X'}_a} = \bm{A\tilde X}_{a} - \bm{B} {\bm{R}_a}^{ - 1}{\bm{B}^{\text{T}}}\left( {\bm{\lambda }\left( f \right) - \bm{\nu }\left( f \right)} \right) \hfill \\
  \bm{\tilde x'}_{da} = \bm{A\tilde x}_{da} + \bm{B}{\bm{R}_d}^{ - 1}{\bm{B}^{\text{T}}}\bm{\nu }\left( f \right) + \bm{B}{\bm{R}_a}^{ - 1}{\bm{B}^{\text{T}}}\left( {\bm{\lambda }\left( f \right) - \bm{\nu }\left( f \right)} \right) \hfill \\ 
  \bm{\lambda }' =  - {\bm{A}^{\text{T}}}\bm{\lambda }\hfill \\
  \bm{\nu }' =  - {\bm{A}^{\text{T}}}\bm{\nu }\hfill \\
\end{gathered}  \right.\label{eq13}
\end{equation}
which can be further expressed:
\begin{equation}
\left[ {\begin{array}{*{20}{c}}
  {\bm{\tilde y}'} \\ 
  {\bm{\Lambda }'} 
\end{array}} \right] = \left[ {\begin{array}{*{20}{c}}
  {{\bm{W}_{11}}}&{{\bm{W}_{12}}} \\ 
  {{\bm{W}_{21}}}&{{\bm{W}_{22}}} 
\end{array}} \right]\left[ {\begin{array}{*{20}{c}}
  {\bm{\tilde y}} \\ 
  \bm{\Lambda } 
\end{array}} \right]
\label{eq14}
\end{equation}
where
\begin{equation}
\begin{gathered}
  {\bm{W}_{11}} = \left[ {\begin{array}{*{20}{c}}
  \bm{A}&{{\bm{O}_{6 \times 6}}} \\ 
  {{\bm{O}_{6 \times 6}}}&\bm{A} 
\end{array}} \right],{\text{ }}{\bm{W}_{12}} = \left[ {\begin{array}{*{20}{c}}
  { - \bm{B}{\bm{R}_a}^{ - 1}{\bm{B}^{\text{T}}}}&{\bm{B}{\bm{R}_a}^{ - 1}{\bm{B}^{\text{T}}}} \\ 
  {\bm{B}{\bm{R}_a}^{ - 1}{\bm{B}^{\text{T}}}}&{\bm{B}\left( {{\bm{R}_d}^{ - 1} - {\bm{R}_a}^{ - 1}} \right){\bm{B}^{\text{T}}}} 
\end{array}} \right], \hfill \\
  {\bm{W}_{21}} = {\bm{O}_{12 \times 12}},{\text{ }}{\bm{W}_{22}} = \left[ {\begin{array}{*{20}{c}}
  { - {\bm{A}^{\text{T}}}}&{{\bm{O}_{6 \times 6}}} \\ 
  {{\bm{O}_{6 \times 6}}}&{ - {\bm{A}^{\text{T}}}} 
\end{array}} \right] \hfill \\ 
\end{gathered}
\label{eq2222}\end{equation}
With the transversality condition presented in Eq. \eqref{eq19}, the solution of Eq. \eqref{eq14} can be expressed as:
\begin{equation}
\left[ {\begin{array}{*{20}{c}}
  {\bm{\tilde y}\left( {{f}} \right)} \\ 
  {\bm{\Lambda }\left( {{f}} \right)} 
\end{array}} \right] = \left[ {\begin{array}{*{20}{c}}
  {{\bm{U }_{11}}\left( {{f},{\text{ }}{f_0}} \right)}&{{\bm{U}_{12}}\left( {{f},{\text{ }}{f_0}} \right)} \\ 
  {{\bm{U}_{21}}\left( {{f},{\text{ }}{f_0}} \right)}&{{\bm{U}_{22}}\left( {{f},{\text{ }}{f_0}} \right)} 
\end{array}} \right]\left[ {\begin{array}{*{20}{c}}
  {\bm{\tilde y}\left( {{f_0}} \right)} \\ 
  {\bm{\Lambda }\left( {{f_0}} \right)} 
\end{array}} \right]\label{eq15}
\end{equation}
where
\begin{equation}
{\bm{U}_{11}}\left( {f,\bm{ }{f_0}} \right) = \left[ {\begin{array}{*{20}{c}}
  {{\bm{\Omega }_{11}}\left( {f,\bm{ }{f_0}} \right)}&{{\bm{O}_{6 \times 6}}} \\ 
  {{\bm{O}_{6 \times 6}}}&{{\bm{\Omega }_{11}}\left( {f,\bm{ }{f_0}} \right)} 
\end{array}} \right] = \left[ {\begin{array}{*{20}{c}}
  {\exp \left( {\int_{{f_0}}^{{f_f}} \bm{A} {\text{d}}f} \right)}&{{\bm{O}_{6 \times 6}}} \\ 
  {{\bm{O}_{6 \times 6}}}&{\exp \left( {\int_{{f_0}}^{{f_f}} \bm{A} {\text{d}}f} \right)} 
\end{array}} \right]
\label{eq2223}\end{equation}
\begin{equation}
{\bm{U}_{22}}\left( {f,\bm{ }{f_0}} \right) = \left[ {\begin{array}{*{20}{c}}
  {{\bm{\Omega }_{22}}\left( {f,\bm{ }{f_0}} \right)}&{{\bm{O}_{6 \times 6}}} \\ 
  {{\bm{O}_{6 \times 6}}}&{{\bm{\Omega }_{22}}\left( {f,\bm{ }{f_0}} \right)} 
\end{array}} \right] = \left[ {\begin{array}{*{20}{c}}
  {\exp \left( {\int_{{f_0}}^{{f}} -\bm{A}^{\text{T}} {\text{d}}f} \right)}&{{\bm{O}_{6 \times 6}}} \\ 
  {{\bm{O}_{6 \times 6}}}&{\exp \left( {\int_{{f_0}}^{{f}} -\bm{A}^{\text{T}} {\text{d}}f} \right)} 
\end{array}} \right]
\label{eq2224}\end{equation}
\begin{equation}
{\bm{U}_{21}}\left( {f,\bm{ }{f_0}} \right) = \bm{O}_{12\times 12}
\label{eq22224}\end{equation}
The matrix ${{\bm{\Omega }_{11}}\left( {{f},{\text{ }}{f_0}} \right)}$ denotes the STM from the initial true anomaly $f_0$ to an epoch $f$ of the uncontrolled TH equations. The equations ${\bm{\tilde X'}_a} = \bm{A\tilde X}_{a} - \bm{B} {\bm{R}_a}^{ - 1}{\bm{B}^{\text{T}}}\left( {\bm{\lambda }\left( f \right) - \bm{\nu }\left( f \right)} \right)$ and $\bm{\tilde x'}_{da} = \bm{A\tilde x}_{da} + \bm{B}{\bm{R}_d}^{ - 1}{\bm{B}^{\text{T}}}\bm{\nu }\left( f \right) + \bm{B}{\bm{R}_a}^{ - 1}{\bm{B}^{\text{T}}}\left( {\bm{\lambda }\left( f \right) - \bm{\nu }\left( f \right)} \right)$ have solutions as follows \cite{li2012fuel}:
\begin{equation}
  {\bm{\tilde X}_a}\left( f \right) = {\bm{\Omega }_{11}}\left( {f,{\text{ }}{f_0}} \right){\bm{\tilde X}_a}\left( {{f_0}} \right) - \int_{{f_0}}^f {{\bm{\Omega }_{11}}\left( {f,{\text{ }}\theta } \right)} \bm{B}{\bm{R}_a}^{ - 1}{\bm{B}^{\text{T}}}\left( {\bm{\lambda }\left( \theta \right) - \bm{\nu }\left( \theta \right)} \right){\text{d}}\theta \label{eq16} 
\end{equation}
\begin{align}\label{eq17}
  {{\bm{\tilde x}}_{da}}\left( f \right) &= {\bm{\Omega }_{11}}\left( {f,{\text{ }}{f_0}} \right){{\bm{\tilde x}}_{da}}\left( {{f_0}} \right) + \int_{{f_0}}^f {{\bm{\Omega }_{11}}\left( {f,{\text{ }}\theta } \right)} \bm{B}{\bm{R}_d}^{ - 1}{\bm{B}^{\text{T}}}\bm{\nu }\left( \theta  \right){\text{d}}\theta  \\ 
   \notag&+ \int_{{f_0}}^f {{\bm{\Omega }_{11}}\left( {f,{\text{ }}\theta } \right)} \bm{B}{\bm{R}_a}^{ - 1}{\bm{B}^{\text{T}}}\left( {\bm{\lambda }\left( \theta  \right) - \bm{\nu }\left( \theta  \right)} \right){\text{d}}\theta  
\end{align}
Since $\bm{\lambda }\left( f \right) = {\bm{\Omega }_{22}}\left( {f,{\text{ }}{f_0}} \right)\bm{\lambda }\left( {{f_0}} \right)$ and $\bm{\nu }\left( f \right) = {\bm{\Omega }_{22}}\left( {f,{\text{ }}{f_0}} \right)\bm{\nu }\left( {{f_0}} \right)$, Eqs. \eqref{eq16}-\eqref{eq17} can be further reconstructed as:
\begin{align}\label{eq1666}
  {\bm{\tilde X}_a}\left( f \right) &= {\bm{\Omega }_{11}}\left( {f,{\text{ }}{f_0}} \right){\bm{\tilde X}_a}\left( {{f_0}} \right) - \int_{{f_0}}^f {{\bm{\Omega }_{11}}\left( {f,{\text{ }}\theta } \right)} \bm{B}{\bm{R}_a}^{ - 1}{\bm{B}^{\text{T}}}{\bm{\Omega }_{22}}\left( {\theta,{\text{ }}f_0 } \right){\text{d}}\theta{\bm{\lambda }\left( f_0 \right)}\\
  \notag&+\int_{{f_0}}^f {{\bm{\Omega }_{11}}\left( {f,{\text{ }}\theta } \right)} \bm{B}{\bm{R}_a}^{ - 1}{\bm{B}^{\text{T}}}{\bm{\Omega }_{22}}\left( {\theta,{\text{ }}f_0 } \right){\text{d}}\theta{\bm{\nu }\left( f_0 \right)} 
\end{align}

\begin{align}\label{eq1667}
  {\bm{\tilde x}_{da}}\left( f \right) &= {\bm{\Omega }_{11}}\left( {f,{\text{ }}{f_0}} \right){\bm{\tilde x}_{da}}\left( {{f_0}} \right) + \int_{{f_0}}^f {{\bm{\Omega }_{11}}\left( {f,{\text{ }}\theta } \right)} \bm{B}{\bm{R}_a}^{ - 1}{\bm{B}^{\text{T}}}{\bm{\Omega }_{22}}\left( {\theta,{\text{ }}f_0 } \right){\text{d}}\theta{\bm{\lambda }\left( f_0 \right)}\\
  \notag&+\int_{{f_0}}^f {{\bm{\Omega }_{11}}\left( {f,{\text{ }}\theta } \right)} \bm{B}\left({\bm{R}_d}^{ - 1}-{\bm{R}_a}^{ - 1}\right){\bm{B}^{\text{T}}}{\bm{\Omega }_{22}}\left( {\theta,{\text{ }}f_0 } \right){\text{d}}\theta{\bm{\nu }\left( f_0 \right)} 
\end{align}
Then, we can obtain:
\begin{equation}
{\bm{U}_{12}}\left(f,\text{ }f_0\right) = \left[ {\begin{array}{*{20}{c}}
  { - {\bm{V}_1}}&{{\bm{V}_1}} \\ 
  {{\bm{V}_1}}&{{\bm{V}_2}} 
\end{array}} \right]\label{eq1668}
\end{equation} 
\begin{equation}
\begin{gathered}
  {\bm{V}_1} = \int_{{f_0}}^f {{\bm{\Omega }_{11}}\left( {f,{\text{ }}\theta } \right)} \bm{B}{\bm{R}_a}^{ - 1}{\bm{B}^{\text{T}}}{\bm{\Omega }_{22}}\left( {\theta ,{\text{ }}{f_{\text{0}}}} \right){\text{d}}\theta  \hfill \\
  {\bm{V}_2} = \int_{{f_0}}^f {{\bm{\Omega }_{11}}\left( {f,{\text{ }}\theta } \right)} \bm{B}\left( {{\bm{R}_d}^{ - 1} - {\bm{R}_a}^{ - 1}} \right){\bm{B}^{\text{T}}}{\bm{\Omega }_{22}}\left( {\theta ,{\text{ }}{f_{\text{0}}}} \right){\text{d}}\theta  \hfill \\ 
\end{gathered}\label{eq1669}
\end{equation}
According to Eq. \eqref{eq15} and the transversality condition, the linear relationship between $\tilde{\bm y}$ and $\bm{\Lambda}$ can be observed:
\begin{equation}
\bm{\Lambda} = \bm{P}\tilde{\bm y}
\label{eq1670}
\end{equation}
The matrix $\bm{P}$ satisfies the matrix Riccati equation \cite{li2024nash}:
\begin{equation}
\bm{P}' = \bm{W}_{22}\bm{P} - \bm{P}\bm{W}_{11} - \bm{P}\bm{W}_{12}\bm{P}
\label{eq16711}
\end{equation}
The transversality condition is presented as follows:
\begin{equation}
\bm{P}\left(f_f\right)  = \left[ {\begin{array}{*{20}{c}}
  {{\bm{S}_a}}&{{\bm{O}_{6 \times 6}}} \\ 
  {{\bm{O}_{6 \times 6}}}&{-{\bm{S}_{da}}} 
\end{array}} \right]
\label{eq16712}
\end{equation}
According to Ref. \cite{li2024nash}, the matrix Riccati equation typically requires numerical backward integration, yielding a numerical LQ strategy for the TAD game. However, conventional numerical solutions face real-time limitations and present difficulties in judging the game outcomes. Therefore, based on the analytical solution of the TH equations developed by Dang \cite{dang2017solutions} and the method proposed by our previous work \cite{fu2025analytical}, we derive the analytical solution of the matrix Riccati equation yielding the analytical strategies and enabling analysis of winning conditions.

\section{Analytical Strategies for the TAD Game}\label{sec3}
In this section, we derive the analytical solution of the matrix Riccati equation and the corresponding analytical strategies for the TAD game.

Combining the transversality condition shown in Eq. \eqref{eq16712} with Eq. \eqref{eq15} and Eq. \eqref{eq1670}, the analytical solution $\bm{P}$ of the Riccati equation at $f$ can be expressed as:
\begin{equation}
\bm{P}\left(f\right)  = \left(\bm{U}_{22}\left(f_f,\text{ }f\right)-\bm{S}\bm{U}_{12}\left(f_f,\text{ }f\right)\right)^{-1}\left(\bm{S}\bm{U}_{11}\left(f_f,\text{ }f\right)-\bm{U}_{21}\left(f_f,\text{ }f\right)\right)
\label{eq1672}
\end{equation}
where
\begin{equation}
\bm{S}  = \left[ {\begin{array}{*{20}{c}}
  {{\bm{S}_a}}&{{\bm{O}_{6 \times 6}}} \\ 
  {{\bm{O}_{6 \times 6}}}&{-{\bm{S}_{da}}} 
\end{array}} \right]
\label{eq1673}
\end{equation}
According to Eq. \eqref{eq1672}, the problem of solving the analytical solution of the matrix Riccati equation can be reduced to obtaining the analytical expressions of the matrices $\bm{U}_{11}$, $\bm{U}_{12}$, and $\bm{U}_{22}$ (note that $\bm{U}_{21}=\bm{O}_{12 \times 12}$). To obtain the analytical expressions, the STM of the uncontrolled TH equations \cite{dang2017solutions} is introduced. ${{\bm{\Omega }_{11}}\left( {{f_2},{\text{ }}{f_1}} \right)}$ shown in Eq. \eqref{eq2223} is expressed as:

\begin{equation}
{\bm{\Omega }_{11}}\left( {{f_2},{\text{ }}{f_1}} \right) = \phi \left( {{f_2}} \right){\phi ^{ - 1}}\left( {{f_1}} \right)\label{eq24}
\end{equation}
The expressions of $\phi \left( f \right)$ and ${\phi ^{ - 1}}\left( f \right)$ can be found in Appendix \ref{secA1}. Therefore, ${{\bm{\Omega }_{22}}\left( {{f_2},{\text{ }}{f_1}} \right)}$ can be expressed as:
\begin{equation}
{\bm{\Omega }_{22}}\left( {{f_2},{\text{ }}{f_1}} \right) = {\left( {{\phi ^{ - 1}}\left( {{f_2}} \right)} \right)^{\text{T}}}{\left( {\phi \left( {{f_1}} \right)} \right)^{\text{T}}}\label{eq25}
\end{equation}
Combining Eq. \eqref{eq24} and Eq. \eqref{eq25} with Eq. \eqref{eq1669}, the matrices $\bm{V}_1$ and $\bm{V}_2$ can be calculated by:
\begin{equation}
{\bm{V }_{1}}\left( {{f_2},{\text{ }}{f_1}} \right) = \frac{1}{{{n^4}}}\left(  \frac{1}{{{r_a}}} \right)\bm{C}_1\label{eq261}
\end{equation}
\begin{equation}
{\bm{V }_{2}}\left( {{f_2},{\text{ }}{f_1}} \right) = \frac{1}{{{n^4}}}\left( {\frac{1}{{{r_d}}} - \frac{1}{{{r_a}}}} \right)\bm{C}_1\label{eq262}
\end{equation}
where
\begin{align}\label{eq27}
\bm{C}_1 &= \phi \left( f_2 \right)\int_{{f_1}}^{f_2} {{\phi ^{ - 1}}\left( \theta  \right)} \left[ {\begin{array}{*{20}{c}}
  \bm{O}&\bm{O} \\ 
  \bm{O}&{\frac{1}{{{\rho ^6}(\theta)}}{\bm{I}_{3 \times 3}}} 
\end{array}} \right]{\phi ^{ - {\text{T}}}}\left( \theta  \right){\text{d}}f  {\phi ^{\text{T}}}\left( {{f_1}} \right) \\
\notag & =\phi \left( f_2 \right)\left({\bm{C}\left(f_2\right)-\bm{C}\left(f_1\right)}\right){\phi ^{\text{T}}}\left( {{f_1}} \right)
\end{align}
The term $\int_{{f_1}}^{f_2} {{\phi ^{ - 1}}\left( \theta  \right)} \left[ {\begin{array}{*{20}{c}}
  \bm{O}&\bm{O} \\ 
  \bm{O}&{\frac{1}{{{\rho ^6}(\theta)}}{\bm{I}_{3 \times 3}}} 
\end{array}} \right]{\phi ^{ - {\text{T}}}}\left( \theta  \right){\text{d}}f$ in Eq. \eqref{eq27} resists direct integration due to the $1/{\rho^i(f)}$ components. Therefore, following the idea of our previous work \cite{fu2025analytical}, we derive the analytical expressions of the aforementioned term, yielding the analytical strategies. Differing from Ref. \cite{fu2025analytical}, we extend the game scenario to the TAD game and further explore the application of these analytical strategies. The eccentric anomaly ($E$) of the reference orbit is introduced to address the aforementioned problem. Geometric relationship between $f$ and $E$ is shown as follows \cite{liu2023algorithms,li2012fuel}:
\begin{equation}
\left\{ \begin{gathered}
  \sin f = \frac{{\sqrt {1 - {e^2}} \sin E}}{{1 - e\cos E}} \hfill \\
  \cos f = \frac{{\cos E - e}}{{1 - e\cos E}} \hfill \\
  L\left( f \right) = nt = \frac{{E - e\sin E}}{{\sqrt {{{\left( {1 - {e^2}} \right)}^3}} }} \hfill \\
  {\text{d}}f = \frac{{\sqrt {1 - {e^2}} }}{{1 - e\cos E}}{\text{d}}E \hfill \\ 
\end{gathered}  \right.\label{eq33}
\end{equation}
The eccentric anomaly can be calculated by the following expressions:
\begin{equation}
E = {\text{atan2}}\left( {\frac{{\sqrt {1 - {e^2}} \sin f}}{{\rho \left( f \right)}},{\text{ }}\frac{{\cos f + e}}{{\rho \left( f \right)}}} \right) + f - {\text{atan2}}\left( {\sin f,{\text{ }}\cos f} \right)\label{eq34}
\end{equation}
This equation is calculated by Matlab® “atan2” command in simulations. Then, Eq. \eqref{eq27} is transformed into \cite{fu2025analytical}:
\begin{align}\label{eqA}
  {\bm{C}_1} &= \phi \left( f{_2} \right)\int_{{E{_1}}}^{E{_2}} {{{\hat \phi }^{ - 1}}\left( \vartheta  \right)} \left[ {\begin{array}{*{20}{c}}
  \bm{O}&\bm{O} \\ 
  \bm{O}&{\frac{1}{{{{\hat \rho }^6}\left( \vartheta  \right)}}{\bm{I}_{3 \times 3}}} 
\end{array}} \right]{{\hat \phi }^{ - {\text{T}}}}\left( \vartheta  \right)\frac{{\sqrt {1 - {e^2}} }}{{1 - e\cos \vartheta }}{\text{d}}\vartheta  \cdot {\phi ^{\text{T}}}\left( {{f{_1}}} \right) \\ 
   \notag &= \phi \left( f{_2} \right)\left( {\bm{C}\left( f{_2} \right) - \bm{C}\left( {{f{_1}}} \right)} \right){\phi ^{\text{T}}}\left( {{f{_1}}} \right) \\ 
   \notag & = \phi \left( f{_2} \right)\left( {\hat{\bm C}\left( E{_2} \right) - \hat{\bm C}\left( {{E{_1}}} \right)} \right){\phi ^{\text{T}}}\left( {{f{_1}}} \right)
\end{align}
where the superscript $\hat {\left(  \cdot  \right)}$ denotes the transformed function from $f$ to $E$. The specific expressions of the matrix $\bm{C}$ can be found in Appendix \ref{secA2}. Notably, the expressions can also be applied to the TAD game along circular reference orbits when the value of $e$ is set to $0$. The expressions are identical to cases along elliptic reference orbits in our previous work \cite{fu2025analytical}. This derivation also confirms that our developed solutions in Ref. \cite{fu2025analytical} can be further extended to multi-player games, such as the TAD game. Combining Eq. \eqref{eq2223}, Eq. \eqref{eq2224}, and Eq. \eqref{eq1668} with Eq. \eqref{eq1672}, we can obtain the analytical solution of the matrix Riccati equation. Once the analytical solution of the matrix Riccati equation is obtained, the analytical Nash-equilibrium strategies can be developed when the state $\tilde{\bm y}$ is known ($\bm{P} = \left[ {\begin{array}{*{20}{c}}
  {{\bm{P}_{11}}}&{{\bm{P}_{12}}} \\ 
  {{\bm{P}_{21}}}&{{\bm{P}_{22}}} 
\end{array}} \right]$, the matrices $\bm{P}_{11}$, $\bm{P}_{12}$, $\bm{P}_{21}$, and $\bm{P}_{22}$ are all $6 \times 6$ matrices):
\begin{equation}
\left\{ \begin{gathered}
   {\bm{u}_a}^* =- {\bm{R}_a}^{ - 1}{\bm{B}^{\text{T}}}\left[\left(\bm{P}_{11}-\bm{P}_{21}\right)\tilde{\bm X}_a+\left(\bm{P}_{12}-\bm{P}_{22}\right)\tilde{\bm x}_{da} \right] \hfill \\
  {\bm{u}_d}^* = {\bm{R}_d}^{ - 1}{\bm{B}^{\text{T}}}\left(\bm{P}_{21}\tilde{\bm X}_a+\bm{P}_{22}\tilde{\bm x}_{da}\right) \hfill \\ 
\end{gathered}  \right.\label{eq3444}
\end{equation}
Since ${\bm{u}_a}^*$ and ${\bm{u}_d}^*$ depend on the state $\tilde{\bm y}$, the strategies are feedback strategies. Notably, to ensure the saddle-point property of the developed analytical strategies, we focus on the TAD game without the constraints of the control inputs.

\noindent\textit{Remark 1:} Compared to the numerical Nash-equilibrium strategies developed by Li et al. \cite{li2024nash}, the developed analytical strategies can effectively improve the computational efficiency.

\noindent\textit{Remark 2:} Based on the similar method proposed by our previous work \cite{fu2025analytical}, we have successfully extended the analytical strategies of orbital pursuit-evasion game to the multi-player game, such as the considered TAD game. To further explore the application of the developed analytical strategies (not only improving the computational efficiency), the analytical winning conditions for the attacker are derived based on the developed strategies.

\noindent\textit{Remark 3:} We focus on the TAD game with a non-maneuvering target \cite{li2024nash} and develop the analytical strategies for it in this paper. The game scenario can be further extended to the TAD game with the target having maneuvering capabilities \cite{li2024orbital} by the similar method developed by this paper and our previous work \cite{fu2025analytical}.

\noindent\textit{Limitation:} Similar to limitations clarified by our previous work \cite{fu2025analytical}, when the eccentricity of the reference orbit is close to 1, the term $\left(1-e^2\right)^j,\text{ }j<0$ could induce numerical instability, which further affects the accuracy of the control inputs. However, for practical applications, the value of eccentricity of the reference orbit is typically low, i.e., $0\leq e \leq 0.8$ \cite{kidder1990use}, which is sufficient for the accuracy of the analytical strategies.

\section{Analytical Winning Conditions for the TAD Game}\label{sec4}
In this section, we present the application of the developed analytical strategies and distinguish our contributions from previous work \cite{fu2025analytical} by deriving winning conditions for the attacker. Due to the three players involved in the game, it proves more complex to analyze the solution of the TAD game than the two-player game \cite{li2024nash}. However, our analytical strategies developed in Section \ref{sec3} enable direct derivation of winning conditions without numerical solutions. In this section, we focus on the case where the attacker and defender have the hovering initial states (i.e., zero initial relative velocity \cite{li2024nash}). The terminal sets of the TAD game are first introduced to define four different game outcomes. Then, the analysis reveals the dependence of the solution on the initial states, leading to analytical winning conditions expressed through the constraints of the initial states.
\subsection{Terminal Set}\label{subsec4.1}
We focus on the TAD game with the fixed game time, i.e., the terminal true anomaly $f_f$ is given \cite{li2024nash} and does not require estimation like Refs. \cite{liao2021research,li2024orbital,fu2025analytical}. To judge the game outcome, we introduce two types of terminal sets as follows \cite{shi2024strategy,li2024nash} (${\bm{\tilde X}}_a=\left[{{\bm{\tilde R}}_a}^\text{T},\text{ }{{\bm{\tilde V}}_a}^\text{T}\right]^\text{T}$ and ${\bm{\tilde x}}_{da}=\left[{{\bm{\tilde r}}_{da}}^\text{T},\text{ }{{\bm{\tilde v}}_{da}}^\text{T}\right]^\text{T}$):
\begin{equation}
\mathcal{T}_1=\left\{\bm{\tilde y}| {\left\| {{{\bm{\tilde R}}_a}} \right\| \leq {R_1}} \right\}
\label{eq44441}
\end{equation}

\begin{equation}
\mathcal{T}_2=\left\{\bm{\tilde y}| {\left\| {{{\bm{\tilde r}}_{da}}} \right\| \leq {R_2}} \right\}
\label{eq44442}
\end{equation}
where $R_1$ denotes the safety distance between the attacker and target, while $R_2$ denotes the safety distance between the attacker and defender. Once the distance between the attacker and target is less than $R_1$, i.e., the state $\bm{\tilde y}$ satisfies $\mathcal{T}_1$, the attacker captures the target successfully; similarly, once the distance between the attacker and defender is less than $R_2$, i.e., the state $\bm{\tilde y}$ satisfies $\mathcal{T}_2$, the defender intercepts the attacker successfully. Therefore, four different outcomes of the considered TAD game are defined:
\begin{enumerate}[label=(\arabic*)]
\item The state $\bm{\tilde y}$ satisfies $\mathcal{T}_1$ at some epoch $f_a \in \left(f_0,\text{ }f_f\right]$ but does not satisfy $\mathcal{T}_2$ at any $f \in \left(f_0,\text{ }f_a\right]$, the attacker is considered as the winner.
\item The state $\bm{\tilde y}$ satisfies $\mathcal{T}_2$ at some epoch $f_d \in \left(f_0,\text{ }f_f\right]$ but does not satisfy $\mathcal{T}_1$ at any epoch $f \in \left(f_0,\text{ }f_d\right]$, the defender is considered as the winner.
\item The state $\bm{\tilde y}$ satisfies both $\mathcal{T}_1$ and $\mathcal{T}_2$ at the same epoch. In this case, both the attacker and defender are not considered as the winner. 
\item The terminal state does not satisfies either $\mathcal{T}_1$ or $\mathcal{T}_2$ at any epoch $f \in \left(f_0,\text{ }f_f\right]$. In this case, we say that both the attacker and defender lose the game.
\end{enumerate}

Based on the aforementioned discussion, we focus on the first case, deriving the analytical sufficient and necessary winning conditions for the attacker with the hovering initial states.
\subsection{Analytical Sufficient And Necessary Winning Conditions for Attacker}\label{subsec4.2}
In this subsection, we present the analytical sufficient and necessary winning conditions for the attacker. Specifically, we focus on the TAD game where the attacker and defender have the hovering initial states. In the considered TAD game with a fixed interval of true anomaly, the dynamical equations of the TAD game Eq. \eqref{eq14} have an analytical solution because the initial costate variables can be analytically obtained by Eq. \eqref{eq1670} and Eq. \eqref{eq1672}:

\begin{equation}
\bm{\Lambda }\left( {{f_0}} \right) = \bm{P}\left(f_0\right)\bm{\tilde y}\left( {{f_0}} \right)\label{new_1}
\end{equation}

\begin{align}\label{new_2}
\bm{\tilde y}\left( f \right) &= {\bm{U}_{11}}\left( {f,\bm{ }{f_0}} \right)\bm{\tilde y}\left( {{f_0}} \right) + {\bm{U}_{12}}\left( {f,\bm{ }{f_0}} \right)\bm{\Lambda }\left( {{f_0}} \right) \\ 
  \notag & = {\bm{U}_{11}}\left( {f,\bm{ }{f_0}} \right)\bm{\tilde y}\left( {{f_0}} \right) + {\bm{U}_{12}}\left( {f,\bm{ }{f_0}} \right)\bm{P}\left( {{f_0}} \right)\bm{\tilde y}\left( {{f_0}} \right) \\ 
  \notag & = \left[ {{\bm{U}_{11}}\left( {f,\bm{ }{f_0}} \right) + {\bm{U}_{12}}\left( {f,\bm{ }{f_0}} \right)\bm{P}\left( {{f_0}} \right)} \right]\bm{\tilde y}\left( {{f_0}} \right) \\ 
  \notag & = \bm{D}\left( {f,{\text{ }}{f_f},{\text{ }}{f_0}} \right)\bm{\tilde y}\left( {{f_0}} \right)
\end{align}
According to Eq. \eqref{new_2}, we can obtain:
\begin{equation}
\left[ {\begin{array}{*{20}{c}}
  {{{\bm{\tilde X}}_a}\left( f \right)} \\ 
  {{{\bm{\tilde x}}_{da}}\left( f \right)} 
\end{array}} \right] = \left[ {\begin{array}{*{20}{c}}
  {{\bm{D}_{11}}\left( {f,\bm{ }{f_f},{\text{ }}{f_0}} \right)}&{{\bm{D}_{12}}\left( {f,\bm{ }{f_f},{\text{ }}{f_0}} \right)} \\ 
  {{\bm{D}_{21}}\left( {f,\bm{ }{f_f},{\text{ }}{f_0}} \right)}&{{\bm{D}_{22}}\left( {f,\bm{ }{f_f},{\text{ }}{f_0}} \right)} 
\end{array}} \right]\left[ {\begin{array}{*{20}{c}}
  {{{\bm{\tilde X}}_a}\left( {{f_0}} \right)} \\ 
  {{{\bm{\tilde x}}_{da}}\left( {{f_0}} \right)} 
\end{array}} \right]
\label{new_3}
\end{equation}
where the matrices ${{\bm{D}_{11}}\left( {f,\bm{ }{f_f},{\text{ }}{f_0}} \right)}$, ${{\bm{D}_{12}}\left( {f,\bm{ }{f_f},{\text{ }}{f_0}} \right)}$, ${{\bm{D}_{21}}\left( {f,\bm{ }{f_f},{\text{ }}{f_0}} \right)}$, and ${{\bm{D}_{22}}\left( {f,\bm{ }{f_f},{\text{ }}{f_0}} \right)}$ are all $6\times 6$ matrices. According to Eq. \eqref{new_3}, the states between the attacker and target can be expressed as:
\begin{equation}
{\bm{\tilde X}_a}\left( {{f}} \right) = {\bm{D}_{11}}\left( {f,\bm{ }{f_f},{\text{ }}{f_0}} \right){\bm{\tilde X}_a}\left( {{f_0}} \right) + {\bm{D}_{12}}\left( {f,\bm{ }{f_f},{\text{ }}{f_0}} \right){\bm{\tilde x}_{da}}\left( {{f_0}} \right)
\label{new_4}
\end{equation}
Therefore, we can conclude that in the considered TAD game using the LQ method, the states between the attacker and target only depend on the initial states of the attacker and defender (i.e., $\bm{\tilde{X}}_a\left(f_0\right)$ and $\bm{\tilde{x}}_{da}\left(f_0\right)$). Equation \eqref{new_4} can be rewritten as:

\begin{align}\label{new_5}
  \left[ {\begin{array}{*{20}{c}}
  {{{\bm{\tilde R}}_a}\left( {{f}} \right)} \\ 
  {{{\bm{\tilde V}}_a}\left( {{f}} \right)} 
\end{array}} \right] &= \left[ {\begin{array}{*{20}{c}}
  {{\bm{D}_{11rr}}\left( {{f},{\text{ }}{f_f},{\text{ }}{f_0}} \right)}&{{\bm{D}_{11rv}}\left( {{f},{\text{ }}{f_f},{\text{ }}{f_0}} \right)} \\ 
  {{\bm{D}_{11vr}}\left( {{f},{\text{ }}{f_f},{\text{ }}{f_0}} \right)}&{{\bm{D}_{11vv}}\left( {{f},{\text{ }}{f_f},{\text{ }}{f_0}} \right)} 
\end{array}} \right]\left[ {\begin{array}{*{20}{c}}
  {{{\bm{\tilde R}}_a}\left( {{f_0}} \right)} \\ 
  {{{\bm{\tilde V}}_a}\left( {{f_0}} \right)} 
\end{array}} \right] \\ 
  \notag & + \left[ {\begin{array}{*{20}{c}}
  {{\bm{D}_{12rr}}\left( {{f},{\text{ }}{f_f},{\text{ }}{f_0}} \right)}&{{\bm{D}_{12rv}}\left( {{f},{\text{ }}{f_f},{\text{ }}{f_0}} \right)} \\ 
  {{\bm{D}_{12vr}}\left( {{f},{\text{ }}{f_f},{\text{ }}{f_0}} \right)}&{{\bm{D}_{12vv}}\left( {{f},{\text{ }}{f_f},{\text{ }}{f_0}} \right)} 
\end{array}} \right]\left[ {\begin{array}{*{20}{c}}
  {{{\bm{\tilde r}}_{da}}\left( {{f_0}} \right)} \\ 
  {{{\bm{\tilde v}}_{da}}\left( {{f_0}} \right)} 
\end{array}} \right] 
\end{align}

\begin{align}\label{new_6}
{{\bm{\tilde R}}_a}\left( {{f}} \right) &= {\bm{D}_{11rr}}\left( {{f},{\text{ }}{f_f},{\text{ }}{f_0}} \right){{\bm{\tilde R}}_a}\left( {{f_0}} \right) + {\bm{D}_{11rv}}\left( {{f},{\text{ }}{f_f},{\text{ }}{f_0}} \right){{\bm{\tilde V}}_a}\left( {{f_0}} \right) \\ 
  \notag & + {\bm{D}_{12rr}}\left( {{f},{\text{ }}{f_f},{\text{ }}{f_0}} \right){{\bm{\tilde r}}_{da}}\left( {{f_0}} \right) + {\bm{D}_{12rv}}\left( {{f},{\text{ }}{f_f},{\text{ }}{f_0}} \right){{\bm{\tilde v}}_{da}}\left( {{f_0}} \right)
\end{align}
Since we focus on the hovering initial states, Eq. \eqref{new_6} can be simplified into:
\begin{equation}
{{\bm{\tilde R}}_a}\left( {{f}} \right) = {\bm{D}_{11rr}}\left( {{f},{\text{ }}{f_f},{\text{ }}{f_0}} \right){{\bm{\tilde R}}_a}\left( {{f_0}} \right) + {\bm{D}_{12rr}}\left( {{f},{\text{ }}{f_f},{\text{ }}{f_0}} \right){{\bm{\tilde r}}_{da}}\left( {{f_0}} \right)
\label{new_7}
\end{equation}
If the matrix ${\bm{D}_{12rr}}\left( {{f},{\text{ }}{f_f},{\text{ }}{f_0}} \right)$ is invertible, Eq. \eqref{new_7} can be expressed as:
\begin{equation}
{\bm{\tilde r}_{da}}\left( {{f_0}} \right) = {\bm{D}_{12rr}}^{ - 1}\left( {f,\bm{ }{f_f},{\text{ }}{f_0}} \right){\bm{\tilde R}_a}\left( {{f}} \right)  - {\bm{D}_{12rr}}^{ - 1}\left( {f,\bm{ }{f_f},{\text{ }}{f_0}} \right){\bm{D}_{11rr}}\left( {f,\bm{ }{f_f},{\text{ }}{f_0}} \right){\bm{\tilde R}_a}\left( {{f_0}} \right)
\label{new_8}
\end{equation}
If the attacker wins the game, the position of the attacker should satisfy $\mathcal{T}_1$ at some epoch $f_a$ satisfying $f_a \in \left(f_0,\text{ }f_f\right]$, i.e., ${\left\| {{{\bm{\tilde R}}_a}\left( {{f_a}} \right)} \right\| \leq {R_1}}$. This illustrates that the position of the attacker at $f_a$ should be constrained inside a sphere:
\begin{equation}
{\bm{\tilde R}_a}\left( {{f_a}} \right) \in \left\{ {\sqrt {{{\tilde x}_a}^2\left( {{f_a}} \right) + {{\tilde y}_a}^2\left( {{f_a}} \right) + {{\tilde z}_a}^2\left( {{f_a}} \right)}  \leq {R_1}} \right\}
\label{new_9}
\end{equation}
Based on the aforementioned discussion, we can obtain that if the position of the attacker satisfies $\mathcal{T}_1$ at $f_a$, the initial position of the defender ${\bm{\tilde R}_d}\left( {{f_0}} \right)={\bm{\tilde R}_a}\left( {{f_0}} \right)+{\bm{\tilde r}_{da}}\left( {{f_0}} \right)$ should be constrained inside an ellipsoid due to the role of the matrix $\bm{D}_{12rr}^{ - 1}\left(f_a,\text{ }f_f,\text{ }f_0\right)$ in stretching and rotation of vector \cite{short2015stretching}. This ellipsoid satisfies the following equation:
\begin{equation}
\left\{ \begin{gathered}
  \bm{\tilde R} = {\left[ {\tilde x,{\text{ }}\tilde y,{\text{ }}\tilde z} \right]^{\text{T}}}|{{\bm{\tilde R}}^{\text{T}}}{\bm{G}_1}\left(f_a\right)\bm{\tilde R} - 2{{\bm{\tilde r}}_1}^{\text{T}}\left(f_a\right){\bm{G}_1}\left(f_a\right)\bm{\tilde R} + {{\bm{\tilde r}}_1}^{\text{T}}\left(f_a\right){\bm{G}_1}\left(f_a\right){{\bm{\tilde r}}_1}\left(f_a\right) - {R_1}^2 = 0 \hfill \\
  {\bm{G}_1}\left(f\right) = {\bm{D}^{\text{T}}}_{12rr}\left( {{f},{\text{ }}{f_f},{\text{ }}{f_0}} \right){\bm{D}_{12rr}}\left( {{f},{\text{ }}{f_f},{\text{ }}{f_0}} \right) \hfill \\
  {{\bm{\tilde r}}_1}\left(f\right) = {{\bm{\tilde R}}_a}\left( {{f_0}} \right) - {\bm{D}^{ - 1}}_{12rr}\left( {{f},{\text{ }}{f_f},{\text{ }}{f_0}} \right){\bm{D}_{11rr}}\left( {{f},{\text{ }}{f_f},{\text{ }}{f_0}} \right){{\bm{\tilde R}}_a}\left( {{f_0}} \right) \hfill \\ 
\end{gathered}  \right.
\label{new_10}
\end{equation}
Notably, to satisfy $\mathcal{T}_1$ at $f_a$, the initial position of the defender should be constrained inside this ellipsoid. Moreover, to ensure that the attacker wins the game, we hope that the distance between the attacker and defender does not satisfy $\mathcal{T}_2$ at any epoch $f\in \left(f_0,\text{ }f_a\right]$. We first derive the initial position of the defender that satisfies $\mathcal{T}_2$. Similar to the aforementioned derivation of the initial position of the defender that satisfies $\mathcal{T}_1$, we have the following relationship:
\begin{equation}
{\bm{\tilde x}_{da}}\left( {{f}} \right) = {\bm{D}_{21}}\left( {f,\bm{ }{f_f},{\text{ }}{f_0}} \right){\bm{\tilde X}_a}\left( {{f_0}} \right) + {\bm{D}_{22}}\left( {f,\bm{ }{f_f},{\text{ }}{f_0}} \right){\bm{\tilde x}_{da}}\left( {{f_0}} \right)
\label{new_11}
\end{equation}

\begin{align}\label{new_12}
  \left[ {\begin{array}{*{20}{c}}
  {{{\bm{\tilde r}}_{da}}\left( {{f}} \right)} \\ 
  {{{\bm{\tilde v}}_{da}}\left( {{f}} \right)} 
\end{array}} \right] &= \left[ {\begin{array}{*{20}{c}}
  {{\bm{D}_{21rr}}\left( {{f},{\text{ }}{f_f},{\text{ }}{f_0}} \right)}&{{\bm{D}_{21rv}}\left( {{f},{\text{ }}{f_f},{\text{ }}{f_0}} \right)} \\ 
  {{\bm{D}_{21vr}}\left( {{f},{\text{ }}{f_f},{\text{ }}{f_0}} \right)}&{{\bm{D}_{21vv}}\left( {{f},{\text{ }}{f_f},{\text{ }}{f_0}} \right)} 
\end{array}} \right]\left[ {\begin{array}{*{20}{c}}
  {{{\bm{\tilde R}}_a}\left( {{f_0}} \right)} \\ 
  {{{\bm{\tilde V}}_a}\left( {{f_0}} \right)} 
\end{array}} \right] \\ 
  \notag & + \left[ {\begin{array}{*{20}{c}}
  {{\bm{D}_{22rr}}\left( {{f},{\text{ }}{f_f},{\text{ }}{f_0}} \right)}&{{\bm{D}_{22rv}}\left( {{f},{\text{ }}{f_f},{\text{ }}{f_0}} \right)} \\ 
  {{\bm{D}_{22vr}}\left( {{f},{\text{ }}{f_f},{\text{ }}{f_0}} \right)}&{{\bm{D}_{22vv}}\left( {{f},{\text{ }}{f_f},{\text{ }}{f_0}} \right)} 
\end{array}} \right]\left[ {\begin{array}{*{20}{c}}
  {{{\bm{\tilde r}}_{da}}\left( {{f_0}} \right)} \\ 
  {{{\bm{\tilde v}}_{da}}\left( {{f_0}} \right)} 
\end{array}} \right] 
\end{align}

\begin{align}\label{new_13}
{{\bm{\tilde r}}_{da}}\left( {{f}} \right) &= {\bm{D}_{21rr}}\left( {{f},{\text{ }}{f_f},{\text{ }}{f_0}} \right){{\bm{\tilde R}}_a}\left( {{f_0}} \right) + {\bm{D}_{21rv}}\left( {{f},{\text{ }}{f_f},{\text{ }}{f_0}} \right){{\bm{\tilde V}}_a}\left( {{f_0}} \right) \\ 
  \notag & + {\bm{D}_{22rr}}\left( {{f},{\text{ }}{f_f},{\text{ }}{f_0}} \right){{\bm{\tilde r}}_{da}}\left( {{f_0}} \right) + {\bm{D}_{22rv}}\left( {{f},{\text{ }}{f_f},{\text{ }}{f_0}} \right){{\bm{\tilde v}}_{da}}\left( {{f_0}} \right)
\end{align}
With the hovering initial states, Eq. \eqref{new_13} can be simplified into:
\begin{equation}
{{\bm{\tilde r}}_{da}}\left( {{f}} \right) = {\bm{D}_{21rr}}\left( {{f},{\text{ }}{f_f},{\text{ }}{f_0}} \right){{\bm{\tilde R}}_a}\left( {{f_0}} \right) + {\bm{D}_{22rr}}\left( {{f},{\text{ }}{f_f},{\text{ }}{f_0}} \right){{\bm{\tilde r}}_{da}}\left( {{f_0}} \right)
\label{new_14}
\end{equation}
Therefore, we can obtain the initial relative position between the attacker and defender when the matrix ${\bm{D}_{22rr}}\left( {{f},{\text{ }}{f_f},{\text{ }}{f_0}} \right)$ is invertible:
\begin{equation}
{\bm{\tilde r}_{da}}\left( {{f_0}} \right) = {\bm{D}_{22rr}}^{ - 1}\left( {f,\bm{ }{f_f},{\text{ }}{f_0}} \right){\bm{\tilde r}_{da}}\left( {{f}} \right)  - {\bm{D}_{22rr}}^{ - 1}\left( {f_f,\bm{ }{f_f},{\text{ }}{f_0}} \right){\bm{D}_{21rr}}\left( {f_f,\bm{ }{f_f},{\text{ }}{f_0}} \right){\bm{\tilde R}_a}\left( {{f_0}} \right)
\label{new_15}
\end{equation}
Similar to the initial position of the defender satisfying $\mathcal{T}_1$, the initial position of the defender satisfying $\mathcal{T}_2$ is also constrained inside an ellipsoid. This ellipsoid can be expressed as:
\begin{equation}
\left\{ \begin{gathered}
  \bm{\tilde R} = {\left[ {\tilde x,{\text{ }}\tilde y,{\text{ }}\tilde z} \right]^{\text{T}}}|{{\bm{\tilde R}}^{\text{T}}}{\bm{G}_2}\left(f\right)\bm{\tilde R} - 2{{\bm{\tilde r}}_2}^{\text{T}}\left(f\right){\bm{G}_2}\left(f\right)\bm{\tilde R} + {{\bm{\tilde r}}_2}^{\text{T}}\left(f\right){\bm{G}_2}\left(f\right){{\bm{\tilde r}}_2}\left(f\right) - {R_2}^2 = 0 \hfill \\
  {\bm{G}_2}\left(f\right) = {\bm{D}^{\text{T}}}_{22rr}\left( {{f},{\text{ }}{f_f},{\text{ }}{f_0}} \right){\bm{D}_{22rr}}\left( {{f},{\text{ }}{f_f},{\text{ }}{f_0}} \right) \hfill \\
  {{\bm{\tilde r}}_2}\left(f\right) = {{\bm{\tilde R}}_a}\left( {{f_0}} \right) - {\bm{D}^{ - 1}}_{22rr}\left( {{f},{\text{ }}{f_f},{\text{ }}{f_0}} \right){\bm{D}_{21rr}}\left( {{f},{\text{ }}{f_f},{\text{ }}{f_0}} \right){{\bm{\tilde R}}_a}\left( {{f_0}} \right) \hfill \\ 
\end{gathered}  \right.
\label{new_16}
\end{equation}
Therefore, if the attacker wins, the initial position of the defender should be located outside the ellipsoid corresponding to $f$ at any epoch $f\in\left(f_0,\text{ }f_a\right]$. Above all, when considering the TAD game with the hovering initial states and assuming that the initial position of the attacker is given, we summarize the winning conditions for the attacker: 

Considering the TAD game satisfying the following assumptions about the initial states:
\begin{enumerate}[label=(\arabic*)]
\item The attacker and defender have the hovering initial states, i.e., ${\bm{\tilde V}_a}\left( {{f_0}} \right) = \bm{0}$ and ${\bm{\tilde v}_{da}}\left( {{f_0}} \right) = \bm{0}$.
\item The initial state $\bm{\tilde y}\left(f_0\right)$ does not satisfy $\mathcal{T}_1$ and $\mathcal{T}_2$, and the initial position of the attacker $\bm{\tilde R}_a\left(f_0\right)$ is given.
\item The attacker and defender adopt the aforementioned Nash-equilibrium strategies without constraints of the control inputs.
\end{enumerate}
The attacker wins the game, if:
\begin{enumerate}[label=(\arabic*)]
\item $\exists f_a\in\left(f_0,\text{ }f_f\right]$, the initial position of the defender $\bm{\tilde R}_d\left(f_a\right)$ satisfies $g_1\left(f_a\right) \leq 0$.
\item $\forall f \in \left(f_0\text{ },f_a\right]$, the initial position of the defender $\bm{\tilde R}_d\left(f\right)$ satisfies $g_2\left(f\right) > 0$.
\end{enumerate}
where
\[g_1\left(f\right)={{\bm{\tilde R}_d\left(f\right)}^{\text{T}}}{\bm{G}_1}\left(f\right)\bm{\tilde R}_d\left(f\right) - 2{{\bm{\tilde r}}_1}^{\text{T}}\left(f\right){\bm{G}_1}\left(f\right)\bm{\tilde R}_d\left(f\right) + {{\bm{\tilde r}}_1}^{\text{T}}\left(f\right){\bm{G}_1}\left(f\right){{\bm{\tilde r}}_1}\left(f\right) - {R_1}^2\]

\[g_2\left(f\right)={{\bm{\tilde R}_d\left(f\right)}^{\text{T}}}{\bm{G}_2}\left(f\right)\bm{\tilde R}_d\left(f\right) - 2{{\bm{\tilde r}}_2}^{\text{T}}\left(f\right){\bm{G}_2}\left(f\right)\bm{\tilde R}_d\left(f\right) + {{\bm{\tilde r}}_2}^{\text{T}}\left(f\right){\bm{G}_2}\left(f\right){{\bm{\tilde r}}_2}\left(f\right) - {R_2}^2\]

Therefore, the set of the initial position of the defender ensuring that the attacker wins the game $\mathcal{S}_a$ can be expressed as:

\begin{equation}
\mathcal{S}_a = \mathop  \cup \limits_{{f_a} \in \left( {{f_0},{\text{ }}{f_f}} \right]} \left( {\mathcal{S}_1\left( {{f_a}} \right)\backslash \left( {\mathop  \cup \limits_{f \in \left( {{f_0},{\text{ }}{f_a}} \right]} \mathcal{S}_2\left( f \right)} \right)} \right)
\label{new_25}
\end{equation}

\begin{equation}
{\mathcal{S}_1}\left(f\right) = \left\{ 
  \left( {\tilde x,{\text{ }}\tilde y,{\text{ }}\tilde z} \right)|{{\bm{\tilde R}}^{\text{T}}}{\bm{G}_1}\left(f\right)\bm{\tilde R} - 2{{\bm{\tilde r}}_1}^{\text{T}}\left(f\right){\bm{G}_1}\left(f\right)\bm{\tilde R} + {{\bm{\tilde r}}_1}^{\text{T}}\left(f\right){\bm{G}_1}\left(f\right){{\bm{\tilde r}}_1}\left(f\right) - {R_1}^2 \leq 0  \right\}
\label{new_261}
\end{equation}

\begin{equation}
{\mathcal{S}_2}\left(f\right) = \left\{ 
  \left( {\tilde x,{\text{ }}\tilde y,{\text{ }}\tilde z} \right)|{{\bm{\tilde R}}^{\text{T}}}{\bm{G}_2}\left(f\right)\bm{\tilde R} - 2{{\bm{\tilde r}}_2}^{\text{T}}\left(f\right){\bm{G}_2}\left(f\right)\bm{\tilde R} + {{\bm{\tilde r}}_2}^{\text{T}}\left(f\right){\bm{G}_2}\left(f\right){{\bm{\tilde r}}_2}\left(f\right) - {R_2}^2 \leq 0  \right\}
\label{new_262}
\end{equation}

\noindent\textit{Remark 4:} Since the aforementioned derivation about the winning conditions are bidirectional (i.e., they hold in both directions, $\Leftrightarrow$), the obtained winning conditions are \textbf{sufficient and necessary conditions}.

\noindent\textit{Remark 5:} The aforementioned derivation focuses on a specific case of the TAD game, i.e., the attacker and defender have the hovering initial states, i.e., their initial velocity is zero \cite{li2024nash}. For more application scenarios (different initial velocity), the analytical winning conditions for the attacker can be further derived based on Eq. \eqref{new_5} and Eq. \eqref{new_12}.

\noindent\textit{Remark 6:} This specific case has been valuably explored by Li et al. \cite{10.1007/978-981-97-3324-8_28,li2024nash}, we have extended their analysis and presented the analytical form of winning conditions. Compared to their analysis based on the numerical strategies, the obtained winning conditions are geometrically intuitive, i.e., when the attacker wins, the initial position of the defender is constrained inside or outside two ellipsoids and their envelopes.

Subsequently, we perform simulations to verify our obtained analytical strategies and winning conditions.

\section{Results and Discussion}\label{sec5}
In this section, simulations of the considered TAD game along the elliptic orbit are performed to verify the effectiveness of the proposed analytical strategies and winning conditions. First, we present an example of the TAD game obtained from the numerical and analytical Nash-equilibrium strategies to verify the effectiveness of the proposed analytical strategies and winning conditions. Then, the effect of the eccentricity on the outcome of the TAD game is analyzed.
\subsection{Effectiveness of the Proposed Strategies}\label{subsec5.1}
In this subsection, an example of the TAD game is first presented to verify the effectiveness of the proposed analytical strategies and winning conditions. Specifically, we focus on the TAD game with the hovering initial states. The parameters used in the simulations are presented in Table \ref{tab1}. When performing simulations using the analytical strategies, we directly obtained $\bm{\tilde X}_a\left(f\right)$ and $\bm{\tilde x}_{da}\left(f\right)$ during the game according to Eq. \eqref{new_2}. In the analytical method, the step-size of $f$ is the same as the step-size to solve the dynamical equation numerically, detailed in the following text. When using the numerical strategies developed by Li et al. \cite{li2024nash}, the numerical integration of the dynamical equation Eq. \eqref{eq13} and the matrix Riccati equation Eq. \eqref{eq16711} is performed. For the numerical method, the fourth Runge-Kutta method is adopted \cite{battin1999an}. The total procedure of the numerical method can be found in Appendix \ref{secA3}. The simulations are performed by an Intel® Core™ i9-13900KF processor with a base frequency of 3.00 GHz, paired with 64 GB of DDR5 RAM.

\begin{table}[!htb]
\caption{Parameter setting for the numerical simulations.}\label{tab1}%
\centering
\renewcommand{\arraystretch}{1.5}
\begin{tabular}{@{}llll@{}}
\hline
Symbol & Value  & Units & Meaning\\
\hline
$\mu$    & $398603$   & $\text{km}{}^3/\text{s}{}^2$  & Earth gravitational constant  \\
$p$    & $10000$   & $\text{km}$  & Semilatus rectum of the reference orbit  \\
$e$    & $0.1$   & --  & Eccentricity of the reference orbit  \\
$f_0$    & $0$   & $\text{rad}$  & Initial true anomaly of the TAD game  \\
$f_f$    & $2\pi$   & $\text{rad}$  & Terminal true anomaly of the TAD game  \\
$h_f$    & $\pi/500$   & $\text{rad}$  & Step-size for solving the dynamical equations  \\
$r_a$    & $5\times10^9$   & --  & Weight of the weighting matrix $\bm{R}_a$  \\
$r_d$    & $3\times10^9$   & --  & Weight of the weighting matrix $\bm{R}_d$  \\
$s_{ar},\text{ }s_{av}$    & $1$   & --  & Weights of the weighting matrix $\bm{S}_a$  \\
$s_{dar},\text{ }s_{dav}$    & $0.001$   & --  & Weights of the weighting matrix $\bm{S}_{da}$  \\
$[\tilde{x}_{a0},\text{ }\tilde{y}_{a0},\text{ }\tilde{z}_{a0}]$    & $[0,\text{ }20,\text{ }0]$   & $\text{km}$  & Initial position of attacker in the LVLH Frame  \\
$[{\tilde x_{a0}}^\prime,\text{ }{\tilde y_{a0}}^\prime,\text{ }{\tilde z_{a0}}^\prime] $    & $[0,\text{ }0,\text{ }0]$   & \text{km/rad}  &Initial velocity of attacker in the LVLH Frame  \\
$[\tilde{x}_{d0},\text{ }\tilde{y}_{d0},\text{ }\tilde{z}_{d0}]$    & $[-2,\text{ }0,\text{ }0]$   & $\text{km}$  & Initial position of defender in the LVLH Frame  \\
$[{\tilde x_{d0}}^\prime,\text{ }{\tilde y_{d0}}^\prime,\text{ }{\tilde z_{d0}}^\prime] $    & $[0,\text{ }0,\text{ }0]$   & \text{km/rad}  &Initial velocity of defender in the LVLH Frame  \\
$R_1$ & $0.01$ & $\text{km}$ & Safe distance between attacker and target \\
$R_2$ & $0.01$ & $\text{km}$ & Safe distance between attacker and defender \\
\hline
\end{tabular}
\end{table}

\begin{figure}[H]
\centering
\includegraphics[width=0.85\textwidth]{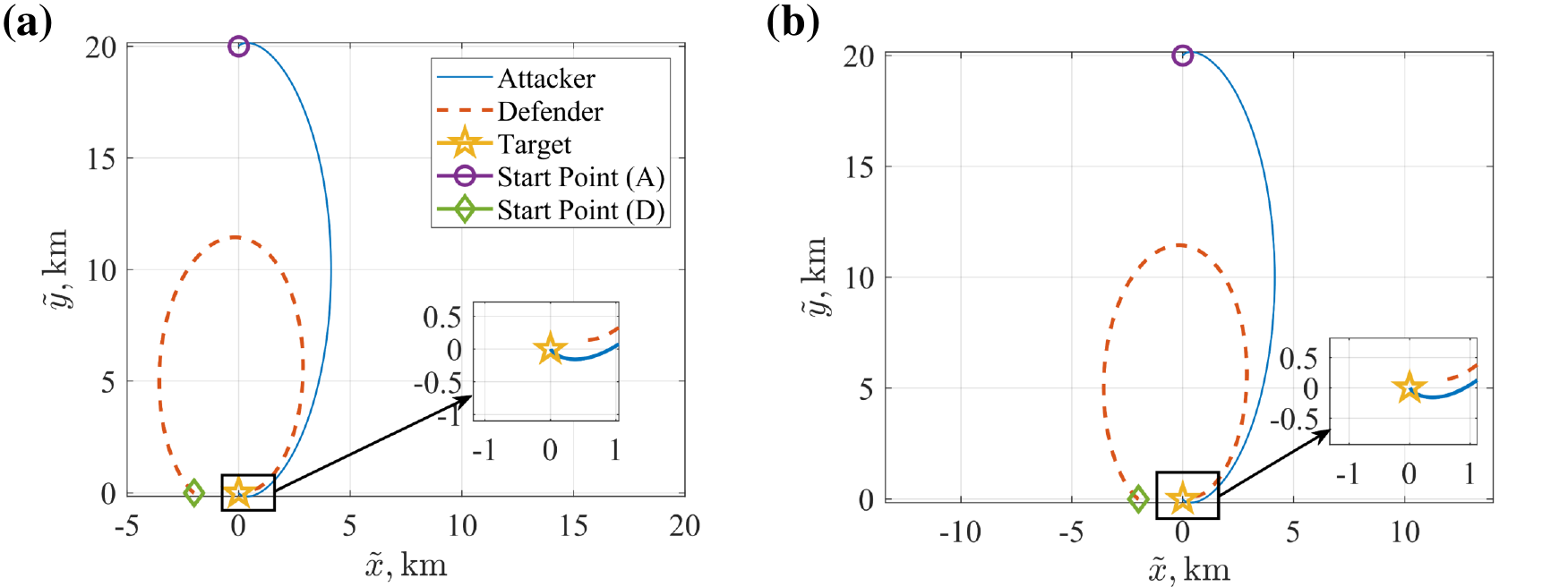}
\caption{The trajectories of the attacker and defender in the LVLH frame. (a) Trajectories obtained from the analytical method; (b) Trajectories obtained from the numerical method.}
\label{fig3}
\end{figure}

\begin{figure}[H]
\centering
\includegraphics[width=0.85\textwidth]{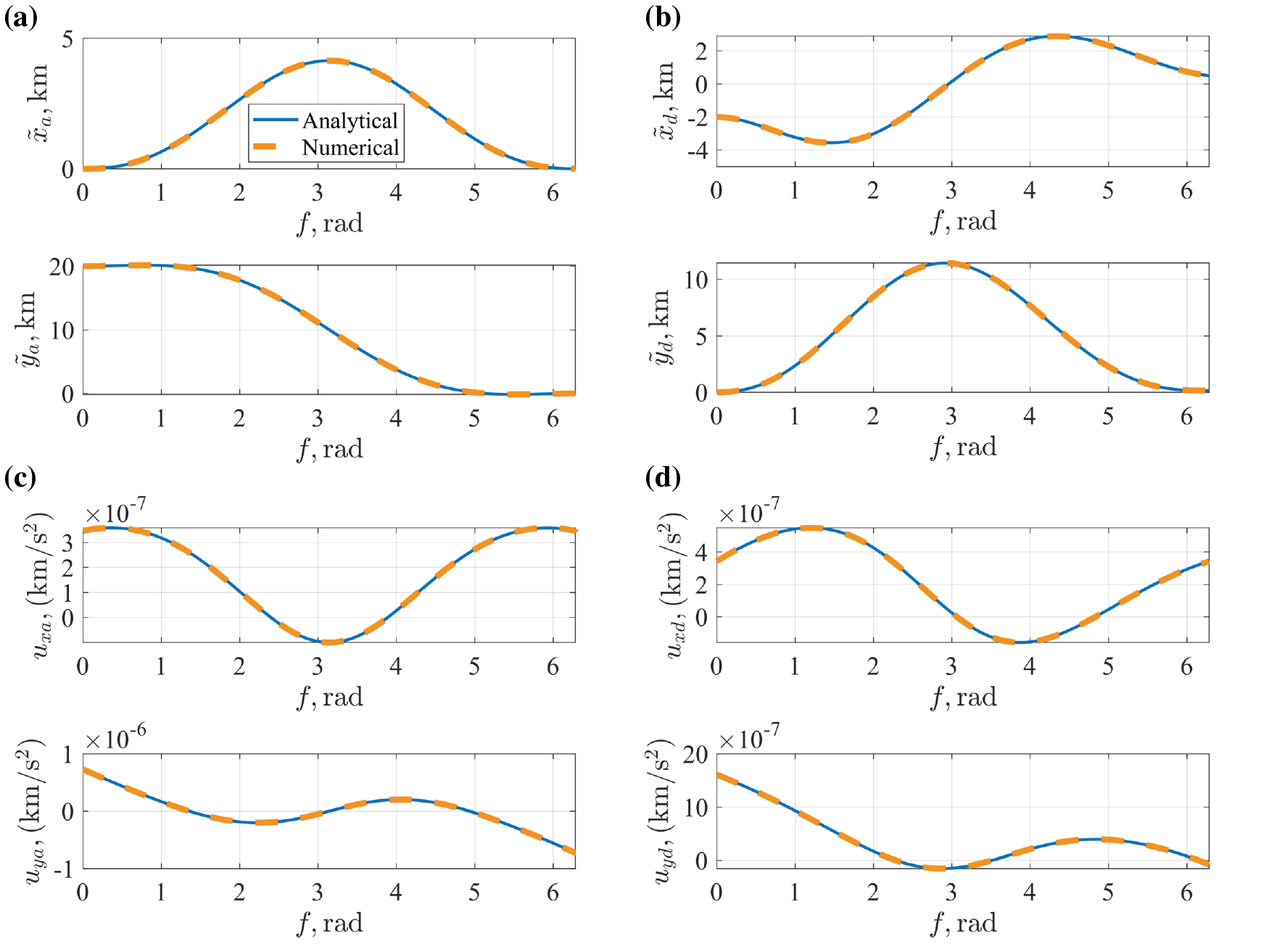}
\caption{The variation in the position and control inputs of attacker and defender during the game. (a) The position of the attacker; (b) The position of the defender; (c) The control inputs of the attacker; (d) The control inputs of the defender.}
\label{fig4}
\end{figure}

The trajectories of the attacker and defender in the LVLH frame obtained from the analytical and numerical methods are shown in Fig. \ref{fig3} (no motion in $\tilde z$-direction). Figure \ref{fig4} presents the variation in position and control inputs of the attacker and defender during the game. From these two figures, we can observe that the analytical solution is in good agreement with the numerical one. The terminal distance between the attacker and the target ($\left|\left|\bm{\tilde R}_a\left(f_f\right)\right|\right|$), terminal distance between the attacker and the defender ($\left|\left|\bm{\tilde r}_{da}\left(f_f\right)\right|\right|$), cost function ($J$), and CPU time obtained from the analytical and numerical methods are summarized in Table \ref{tab2}. The data shown in Table \ref{tab2} confirm the good consistency between the analytical and numerical methods, as the relative error of 0.025$\%$ in $\left|\left|\bm{\tilde R}_a\left(f_f\right)\right|\right|$, 0.008$\%$ in $\left|\left|\bm{\tilde r}_{da}\left(f_f\right)\right|\right|$ and 0.004$\%$ in $J$ is observed. Moreover, the analytical method saves CPU time over 99.9$\%$ than the numerical one, which can be considered as an advantage of the proposed analytical strategies.

\begin{table}[!htb]
\caption{Comparison between analytical and numerical methods.}\label{tab2}%
\centering
\renewcommand{\arraystretch}{1.5}
\begin{tabular}{@{}lllll@{}}
\hline
Method & $\left|\left|\bm{\tilde R}_a\left(f_f\right)\right|\right|,\text{ km}$  & $\left|\left|\bm{\tilde r}_{da}\left(f_f\right)\right|\right|,\text{ km}$ & $J$ & CPU Time, s \\
\hline
Analytical   &  $3.2018 \times 10^{-3}$  & $0.50914$  & $-2.4361\times 10^{-3}$ & $0.10178$  \\
Numerical    &  $3.2010 \times 10^{-3}$  & $0.50910$  & $-2.4362\times 10^{-3}$ & $830.92$  \\
\hline
\end{tabular}
\end{table}

\begin{figure}[H]
\centering
\includegraphics[width=0.95\textwidth]{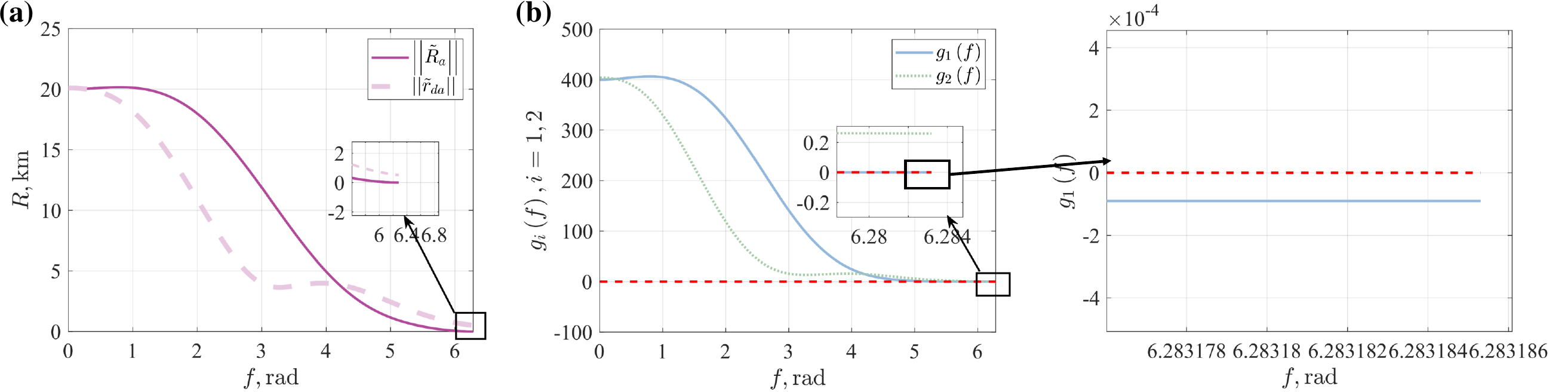}
\caption{The distance of the attacker and target, the distance between attacker and defender, $g_1\left(f\right)$, and $g_2\left(f\right)$ during the game.}
\label{fig5}
\end{figure}

\begin{figure}[H]
\centering
\includegraphics[width=0.75\textwidth]{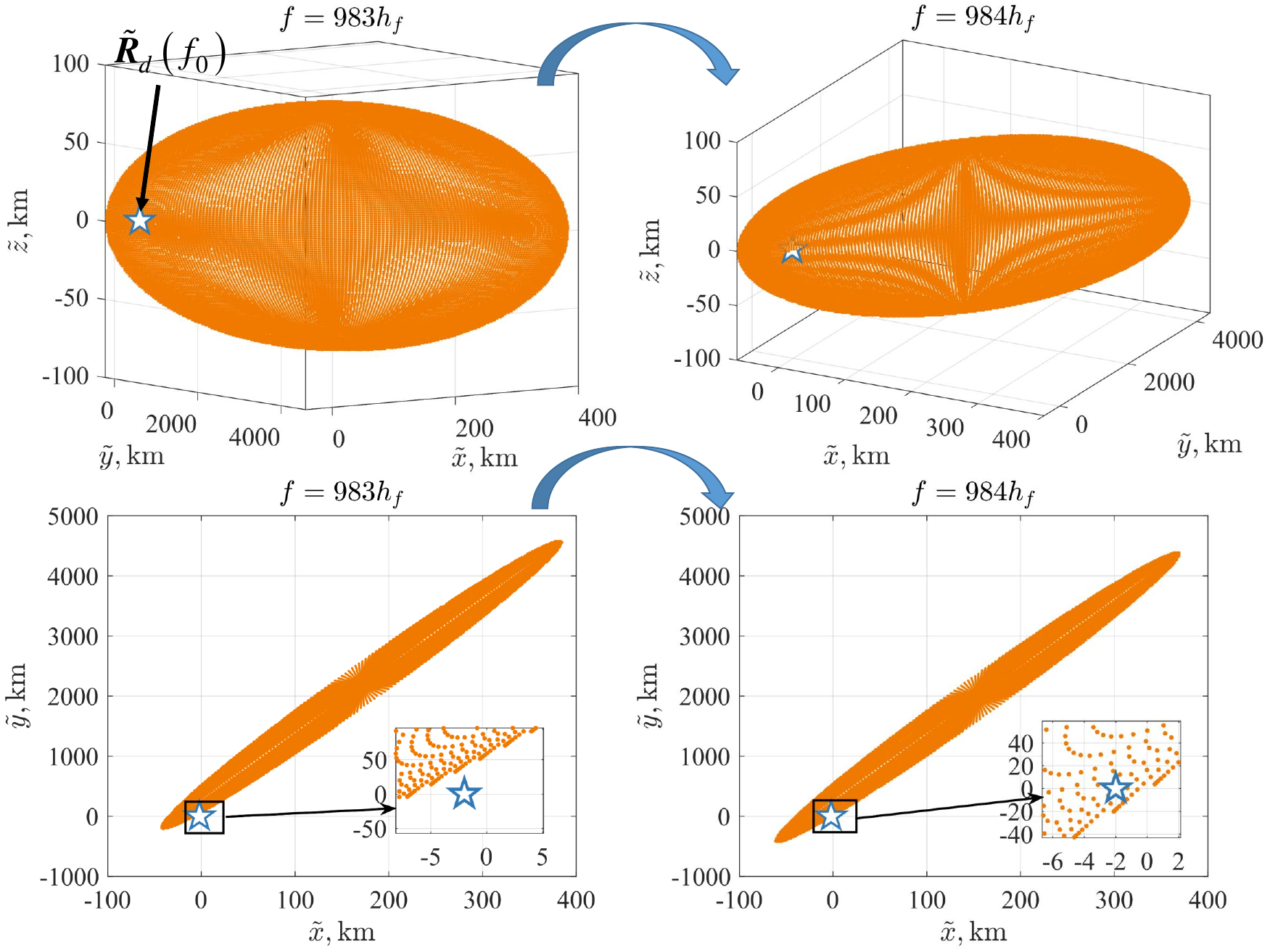}
\caption{The geometric interpretation of the proposed winning conditions.}
\label{fig6}
\end{figure}

To further analyze the game outcome, we present the distance between the attacker and target, the distance between the attacker and defender, $g_1\left(f\right)$, and $g_2\left(f\right)$ during the game in Fig. \ref{fig5} (for the curves of $g_1\left(f\right)$ and $g_2\left(f\right)$, $f\in\left(f_0,\text{ }f_f\right]$). From this figure, we can observe that the curve of $g_1\left(f\right)$ has a zero point, while $g_2\left(f\right)>0$ holds for all time during the game. Therefore, the attacker wins the game. The distance shown in Fig. \ref{fig5} (a) confirms this result. In Fig. \ref{fig5} (b), $g_1\left(f\right)<0$ holds after one specific epoch. Geometrically, there is an epoch $f_{an}$ when $\bm{\tilde R}_d\left(f_0\right)$ is located outside the ellipsoid corresponding to $f_{an}$ defined in Eq. \eqref{new_10}, while at $f_{an}+h_f$, $\bm{\tilde R}_d\left(f_0\right)$ “enters” the ellipsoid corresponding to $f_{an}+h_f$, as shown in Fig. \ref{fig6}. In this game, $f_{an}=983h_f$. Numerically, we treat $f_{an}+h_f$ as $f_a$. These results verify the effectiveness of the proposed winning conditions for the attacker.

Furthermore, we analyze the effects of the eccentricity on the outcome of the considered TAD game.
\subsection{Effects of Eccentricity}\label{subsec5.2}
\begin{figure}[H]
\centering
\includegraphics[width=0.85\textwidth]{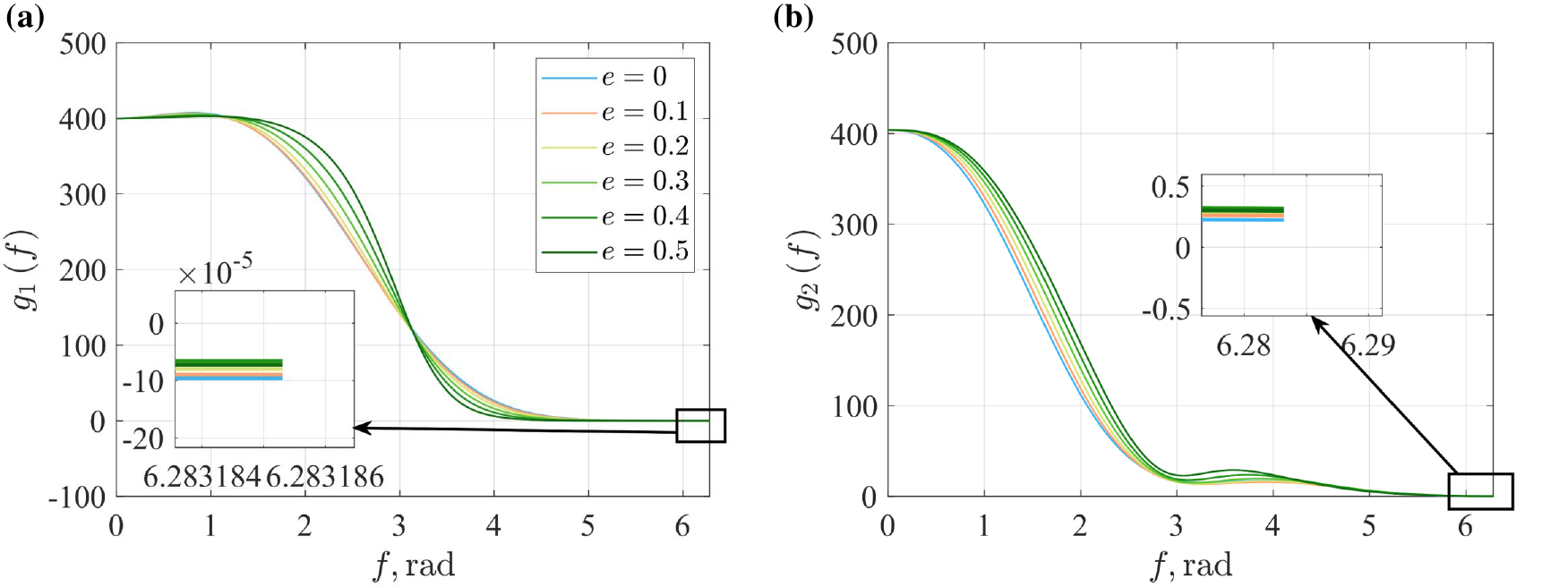}
\caption{The minimum value of $g_1\left(f\right)$ and $g_2\left(f\right)$ during the game under different values of eccentricity. (a) The minimum value of $g_1\left(f\right)$; (b) The minimum value of $g_2\left(f\right)$.}
\label{fig7}
\end{figure}

As the effects of the weights of the weighting matrix on the game outcomes are comprehensively analyzed by Li et al. \cite{li2024nash}, in this paper, we focus on the effects of eccentricity on the game outcomes due to our extension of their model from the CW relative dynamics to the TH relative dynamics. With the same conditions shown in Table \ref{tab1} except for the eccentricity, we select the eccentricity within $\left[0,\text{ }0.5\right]$ with a step-size of 0.1, and present the variation of $g_1\left(f\right)$ and $g_2\left(f\right)$ during the game in Fig. \ref{fig7}. From this figure, it is observed that all the curves of $g_1\left(f\right)$ have zero points, while $g_2\left(f\right)>0$ for all the time. It can be concluded that the eccentricity does not affect the outcome that the attacker wins the game in the considered cases. 

\begin{figure}[H]
\centering
\includegraphics[width=0.4\textwidth]{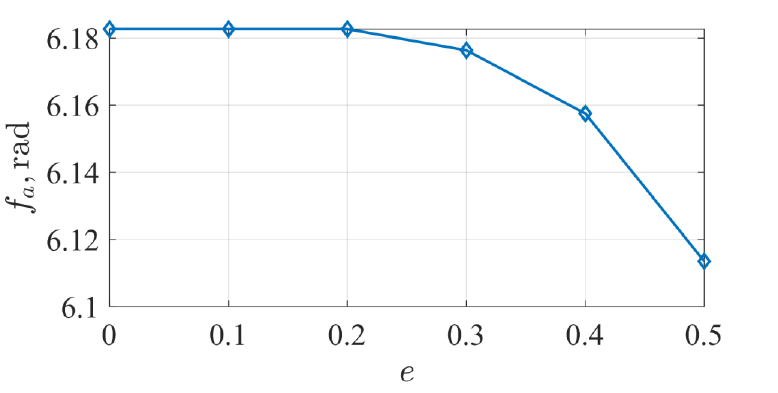}
\caption{The variation of $f_a$ with respect to the eccentricity.}
\label{fig8}
\end{figure}

As $e$ does not affect the game outcome, we further explore the effects of $e$ on the characterization of the game, such as $f_a$. The variation of $f_a$ with respect to eccentricity is shown in Fig. \ref{fig8}. As the value of $e$ increases ($0.2\leq e \leq 0.5$), a gradual decrease in $f_a$ is observed (the case where the value of $f_a$ is not observed to decrease in $e \in \left[0,\text{ }0.2\right]$ probably due to the accuracy of the step-size $h_f$).

\section{Conclusion}\label{sec6}
We focus on the linear quadratic Target-Attacker-Defender (TAD) game with a non-maneuvering target along the elliptic orbits in this paper and devote to providing an analytical framework for it. First, we have derived the analytical Nash-equilibrium strategies for the TAD game, yielding an analytical solution of the game. Based on the proposed analytical strategies and analytical solution of the game, we further obtained the analytical sufficient and necessary winning conditions for the attacker in the TAD game with the hovering initial states. To ensure that the attacker wins the game, we point out that the initial position of the defender should be constrained in a specific set, which can be analytically expressed and interpreted from the geometry through our proposed method. The simulation results verify the effectiveness of the proposed analytical strategies and winning conditions, and the comparison between the results obtained from the analytical and numerical methods further strengthens the advantage of the analytical strategies. The effects of eccentricity on the game outcomes are analyzed. In the considered cases, the eccentricity does not affect the qualitative outcome of the game but affects the epoch when the attacker captures the target. This work provides the analytical framework for orbital Target-Attacker-Defender games, establishing a geometric link between the initial states of players and game outcomes.

\section*{Appendix A: Components of the Matrix $\phi \left( f \right)$ in Elliptic Reference Orbits}\label{secA1}

Components of the Matrix $\phi \left( f \right)$ can be expressed as follows \cite{dang2017solutions}:
\begin{equation}
\phi \left( f \right) = \left[ {\begin{array}{*{20}{c}}
  {\begin{array}{*{20}{c}}
  {{\varphi _1}}&{{\varphi _2}}&{{\varphi _3}}&0&0&0 \\ 
  { - 2S\left( {{\varphi _1}} \right)}&{ - 2S\left( {{\varphi _2}} \right)}&{ - S\left( {2{\varphi _3} + 1} \right)}&1&0&0 \\ 
  0&0&0& 0&{\cos f}&{\sin f} \\
  {{\varphi _1}^\prime }&{{\varphi _2}^\prime }&{{\varphi _3}^\prime }&0&0&0 \\
   { - 2{\varphi _1}}&{ - 2{\varphi _2}}&{ - 2{\varphi _3} - 1}&0&0&0 \\
   0&0&0&0&{ - \sin f}&{\cos f} 
\end{array}} 
\end{array}} \right]\label{eq28}
\end{equation}

\begin{equation}
{\phi ^{ - 1}}\left( f \right) = \left[ {\begin{array}{*{20}{c}}
  {\begin{array}{*{20}{c}}
  {4S\left( {{\varphi _2}} \right) + {\varphi _2}^\prime }&0&0&{ - {\varphi _2}}&{2S\left( {{\varphi _2}} \right)}&0 \\ 
  { - 4S\left( {{\varphi _1}} \right) - {\varphi _1}^\prime }&0&0&{{\varphi _1}}&{ - 2S\left( {{\varphi _1}} \right)}&0 \\ 
  { - 2}&0&0 &0&{ - 1}&0 \\
   { - 2S\left( {2{\varphi _3} + 1} \right) - {\varphi _3}^\prime }&1&0&{{\varphi _3}}&{ - S\left( {2{\varphi _3} + 1} \right)}&0 \\
   0&0&{\cos f}&0&0&{ - \sin f} \\
   0&0&{\sin f}&0&0&{\cos f} 
\end{array}} 
\end{array}} \right]\label{eq29}
\end{equation}

\begin{equation}
\left\{ \begin{gathered}
  {\varphi _1} = \rho \left( f \right)\sin f \hfill \\
  {\varphi _1}^\prime  = \rho \left( f \right)\cos f - e{\sin ^2}f \hfill \\
  S\left( {{\varphi _1}} \right) =  - \cos f - \frac{1}{2}e{\cos ^2}f \hfill \\ 
\end{gathered}  \right.\label{eq30}
\end{equation}
For elliptic reference orbits $\left( {0 \leq e < 1} \right)$,
\begin{equation}
\left\{ \begin{gathered}
  {\varphi _2} = \frac{{e{\varphi _1}}}{{1 - {e^2}}}\left( {D\left( f \right) - 3eL\left( f \right)} \right) - \frac{{\cos f}}{{\rho \left( f \right)}} \hfill \\
  {\varphi _3} =  - \frac{{{\varphi _1}}}{{1 - {e^2}}}\left( {D\left( f \right) - 3eL\left( f \right)} \right) - \frac{{{{\cos }^2}f}}{{\rho \left( f \right)}} - {\cos ^2}f \hfill \\
  {\varphi _2}^\prime  = \frac{{e{\varphi _1}^\prime }}{{1 - {e^2}}}\left( {D\left( f \right) - 3eL\left( f \right)} \right) + \frac{{e\sin f\cos f}}{{{\rho ^2}\left( f \right)}} + \frac{{\sin f}}{{\rho \left( f \right)}} \hfill \\
  {\varphi _3}^\prime  = 2\left( {{\varphi _1}^\prime S\left( {{\varphi _2}} \right) - {\varphi _2}^\prime S\left( {{\varphi _1}} \right)} \right) \hfill \\
  S\left( {{\varphi _2}} \right) =  - \frac{{{\rho ^2}\left( f \right)\left( {D\left( f \right) - 3eL\left( f \right)} \right)}}{{2\left( {1 - {e^2}} \right)}} \hfill \\
  S\left( {2{\varphi _3} + 1} \right) = \frac{{e\sin f\left( {2 + e\cos f} \right)}}{{1 - {e^2}}} - \frac{{3{\rho ^2}\left( f \right)L\left( f \right)}}{{1 - {e^2}}} \hfill \\
  D\left( f \right) = \frac{{\sin f\left( {2 + e\cos f} \right)}}{{{\rho ^2}\left( f \right)}} \hfill \\
  L\left( f \right) = \int {\frac{1}{{{\rho ^2}\left( f \right)}}} {\text{d}}f = nt \hfill \\ 
\end{gathered}  \right.\label{eq31}
\end{equation}

\section*{Appendix B: Components of the Matrix $\bm{C}$ in Elliptic Reference Orbits}\label{secA2}

Since the matrix $\bm{C}$ is a symmetric matrix, only the upper triangular components are presented \cite{fu2025analytical}. The subscript $\bm{F}_{ij}$ denotes the element of the \textit{i}-th row and \textit{j}-th column of the matrix $\bm{F}$. These expressions can also be found in Appendix B of Ref. \cite{fu2025analytical}. Notably, in the origin equations in our previous work, there is a clerical error in Eq. (55), the term $\cos 3E$ should be replaced with ${\cos^3} E$.

\begin{equation}
 {\bm{C}_{11}}\left( f \right) =  - \frac{{3 }}{{2{{\left( {1 - {e^2}} \right)}^{\frac{15}{2}}}}}\left( \begin{gathered}
   - \left( {\left( { - \frac{2}{{15}}{e^4} - \frac{2}{{15}}{e^2} + \frac{2}{{15}}} \right)\sin E + {e^3}E} \right){e^3}{\cos ^4}E \hfill \\
   - \frac{8}{3}\left( {\left( { - \frac{5}{{32}}{e^4} - \frac{1}{{16}}{e^2} - \frac{5}{{16}}} \right)\sin E + e\left( {1 + {e^2}} \right)E} \right){e^2}{\cos ^3}E \hfill \\
   - 2e\left( {\left( {\frac{7}{9} - \frac{4}{{45}}{e^6} + \left( {{E^2} - \frac{{44}}{{45}}} \right){e^4} - \frac{7}{{15}}{e^2}} \right)\sin E - \frac{{13}}{2}{e^3}E - 4eE} \right){\cos ^2}E \hfill \\
   + \left( {\left( {1 - \frac{{49}}{{24}}{e^6} + \left( {3{E^2} - \frac{{187}}{{12}}} \right){e^4} - \frac{{85}}{{12}}{e^2}} \right)\sin E - 8e\left( {1 + {e^4}} \right)E} \right)\cos E \hfill \\
   + \left( {\frac{{16}}{{45}}{e^7} + \left( { - 4{E^2} + \frac{{1016}}{{45}}} \right){e^5} + \left( {\frac{{126}}{5} + 6{E^2}} \right){e^3} + \frac{{128}}{9}e} \right)\sin E \hfill \\
   + \left( { - \frac{5}{3} - \frac{{49}}{{24}}{e^6} + \left( {{E^2} - \frac{{251}}{{12}}} \right){e^4} + \left( { - \frac{{63}}{4} - 2{E^2}} \right){e^2}} \right)E \hfill \\ 
\end{gathered}  \right)
\label{eqA1}
\end{equation}

\begin{equation}
{\bm{C}_{12}}\left( f \right) = \frac{1}{{60{{\left( {1 - {e^2}} \right)}^6}}}\left( \begin{gathered}
   - 12{e^3}{\cos ^5}E + \left( { - 15{e^6} + 45{e^4} + 75{e^2}} \right){\cos ^4}E + \left( { - 40{e^5} - 100{e^3} - 140e} \right){\cos ^3}E \hfill \\
   + \left( { - 60{e^3}\left( {{e^2} - 2} \right)E\sin E + 150{e^4} + 150{e^2} + 90} \right){\cos ^2}E \hfill \\
   + \left( {\left( {90{e^4}E - 360{e^2}E} \right)\sin E - 120{e^5} + 240{e^3} + 300e} \right)\cos E \hfill \\
   + 45\left( {\left( { - {e^3} + 4e} \right){{\sin }^2}E - \frac{8}{3}\left( {{e^4} - \frac{7}{2}{e^2} - 3} \right)E\sin E + e\left( {{e^2} - 6} \right){E^2}} \right)e \hfill \\ 
\end{gathered}  \right)
\label{eqA2}
\end{equation}

\begin{equation}
{\bm{C}_{13}}\left( f \right) = \frac{1}{{60{{\left( {1 - {e^2}} \right)}^6}}}\left( \begin{gathered}
   - 12{e^4}{\cos ^5}E + \left( {30{e^5} + 75{e^3}} \right){\cos ^4}E + \left( { - 100{e^4} - 180{e^2}} \right){\cos ^3}E \hfill \\
   + \left( {60{e^4}E\sin E + 180{e^3} + 210e} \right){\cos ^2}E \hfill \\
   + \left( { - 270{e^3}E\sin E + 120{e^4} + 420{e^2} - 120} \right)\cos E + 135{e^3}{\sin ^2}E \hfill \\
   + \left( {120{e^4}E + 540{e^2}E} \right)\sin E - 135{e^3}{E^2} - 90e{E^2} \hfill \\ 
\end{gathered}  \right)
\label{eqA3}
\end{equation}

\begin{equation}
{\bm{C}_{14}}\left( f \right) =  - \frac{3}{{2{{\left( {1 - {e^2}} \right)}^{\frac{{15}}{2}}}}}\left( \begin{gathered}
  \left( {\left( { - \frac{2}{{15}}{e^4} + \frac{4}{{15}}{e^6}} \right)\sin E - {e^5}E} \right){\cos ^4}E \hfill \\
   + \frac{4}{3}{e^2}{\cos ^3}E\left( {\left( { - \frac{1}{4}{e^5} + \frac{3}{{16}}{e^3} + \frac{9}{8}e} \right)\sin E + E\left( {{e^4} - 4{e^2} - 1} \right)} \right) \hfill \\
   - 2\left( \begin{gathered}
  \left( {\frac{{{e^7}}}{9} - \frac{8}{{45}}{e^5} + {e^3}\left( {{E^2} - \frac{{86}}{{45}}} \right) + \frac{{11}}{9}e} \right)\sin E \hfill \\
   + 2\left( {{e^4} - \frac{{21}}{4}{e^2} - 1} \right)E \hfill \\ 
\end{gathered}  \right)e{\cos ^2}E \hfill \\
   + \left( {\left( {\frac{{11}}{6}{e^7} - \frac{{13}}{8}{e^5} + \left( {3{E^2} - \frac{{275}}{{12}}} \right){e^3} - e} \right)\sin E - 4E{{\left( {1 + {e^2}} \right)}^2}} \right)\cos E \hfill \\
   + \left( {\frac{{16}}{3} - \frac{4}{9}{e^8} - \frac{{298}}{{45}}{e^6} + \left( {\frac{{1304}}{{45}} - 4{E^2}} \right){e^4} + \left( {6{E^2} + \frac{{316}}{9}} \right){e^2}} \right)\sin E \hfill \\
   + eE\left( {\frac{{11}}{6}{e^6} + \frac{{41}}{{24}}{e^4} + \left( {{E^2} - \frac{{129}}{4}} \right){e^2} - 2{E^2} - \frac{{35}}{3}} \right) \hfill \\ 
\end{gathered}  \right)
\label{eqA4}
\end{equation}

\begin{equation}
{\bm{C}_{15}}\left( f \right) =0 \text{ }\text{ }\text{ }\text{ }{\bm{C}_{16}}\left( f \right) =0
\label{eqA5}
\end{equation}

\begin{equation}
{\bm{C}_{22}}\left( f \right) = \frac{{1 }}{{2{{\left( {1 - {e^2}} \right)}^{\frac{11}{2}}}}}\left( \begin{gathered}
  \left( \begin{gathered}
   - \frac{2}{5}{e^3}{\cos ^4}E + \frac{5}{2}{e^2}{\cos ^3}E - \frac{2}{3}e\left( {{e^6} - 5{e^4} + \frac{{39}}{5}{e^2} + 7} \right){\cos ^2}E \hfill \\
   + \left( {{e^6} - 9{e^4} + \frac{{75}}{4}{e^2} + 3} \right)\cos E - \frac{4}{3}{e^7} + \frac{{26}}{3}{e^5} - \frac{{22}}{5}{e^3} - \frac{{82}}{3}e \hfill \\ 
\end{gathered}  \right)\sin E \hfill \\
   + E\left( {{e^6} - 11{e^4} + \frac{{83}}{4}{e^2} + 5} \right) \hfill \\ 
\end{gathered}  \right)\label{eqA6}
\end{equation}

\begin{equation}
{\bm{C}_{23}}\left( f \right) =  - \frac{{3 }}{{2{{\left( {1 - {e^2}} \right)}^{\frac{11}{2}}}}}\left( {\left( \begin{gathered}
  \frac{2}{{15}}{e^4}{\cos ^4}E - \frac{5}{6}{e^3}{\cos ^3}E \hfill \\
   + \left( { - \frac{2}{9}{e^6} + \frac{{28}}{{45}}{e^4} + 2{e^2}} \right){\cos ^2}E \hfill \\
   + \left( {{e^5} - \frac{{13}}{4}{e^3} - \frac{7}{3}e} \right)\cos E \hfill \\
   - \frac{4}{9}{e^6} - \frac{{34}}{{45}}{e^4} + 8{e^2} + \frac{4}{3} \hfill \\ 
\end{gathered}  \right)\sin E + eE\left( {{e^4} - \frac{{31}}{{12}}{e^2} - \frac{{11}}{3}} \right)} \right)\label{eqA7}
\end{equation}

\begin{equation}
{\bm{C}_{24}}\left( f \right) = \frac{1}{{60{{\left( {1 - {e^2}} \right)}^6}}}\left( \begin{gathered}
   - 12{e^4}{\cos ^5}E + \left( { - 30{e^5} + 135{e^3}} \right){\cos ^4}E + \left( {40{e^6} - 100{e^4} - 220{e^2}} \right){\cos ^3}E \hfill \\
   + \left( { - 60{e^2}\left( {{e^2} - 2} \right)E\sin E - 120{e^5} + 465{e^3} - 90e} \right){\cos ^2}E \hfill \\
   + \left( {\left( {90{e^3}E - 360eE} \right)\sin E - 60{e^2} + 480} \right)\cos E - 120\left( {{e^4} - \frac{7}{2}{e^2} - 3} \right)E\sin E \hfill \\
   + 45\left( {\left( {{E^2} - 1} \right){e^2} - 6{E^2} + 4} \right)e \hfill \\ 
\end{gathered}  \right)\label{eqA8}
\end{equation}

\begin{equation}
{\bm{C}_{25}}\left( f \right) =0 \text{ }\text{ }\text{ }\text{ }{\bm{C}_{26}}\left( f \right) =0
\label{eqA9}
\end{equation}

\begin{equation}
{\bm{C}_{33}}\left( f \right) = \frac{1}{{120{{\left( {1 - {e^2}} \right)}^{\frac{{11}}{2}}}}}\left( \begin{gathered}
  \left( \begin{gathered}
   - 24{e^5}{\cos ^4}E + 150{e^4}{\cos ^3}E + \left( { - 32{e^5} - 400{e^3}} \right){\cos ^2}E \hfill \\
   + \left( {225{e^4} + 600{e^2}} \right)\cos E - 64{e^5} - 800{e^3} - 600e \hfill \\ 
\end{gathered}  \right)\sin E \hfill \\
   + 225{e^4}E + 600{e^2}E + 120E \hfill \\ 
\end{gathered}  \right)\label{eqA10}
\end{equation}

\begin{equation}
{\bm{C}_{34}}\left( f \right) = \frac{1}{{60{{\left( {1 - {e^2}} \right)}^6}}}\left( \begin{gathered}
   - 12{e^5}{\cos ^5}E + \left( { - 15{e^6} + 120{e^4}} \right){\cos ^4}E + \left( {60{e^5} - 340{e^3}} \right){\cos ^3}E \hfill \\
   + \left( {60{e^3}E\sin E - 90{e^4} + 345{e^2}} \right){\cos ^2}E + \left( { - 270{e^2}E{{\sin }^2}E + 180{e^3} + 240e} \right)\cos E \hfill \\
   + \left( {120{e^3}E + 540eE} \right)\sin E + \left( { - 135{E^2} + 135} \right){e^2} - 90{E^2} \hfill \\ 
\end{gathered}  \right)\label{eqA11}
\end{equation}

\begin{equation}
{\bm{C}_{35}}\left( f \right) =0 \text{ }\text{ }\text{ }\text{ }{\bm{C}_{36}}\left( f \right) =0
\label{eqA12}
\end{equation}

\begin{equation}
{\bm{C}_{44}}\left( f \right) = \frac{{3eE}}{{{{\left( {1 - {e^2}} \right)}^{\frac{{15}}{2}}}}}\left( \begin{gathered}
  \left( \begin{gathered}
   - \frac{1}{{15}}{e^6}{\cos ^4}E + \left( { - \frac{1}{6}{e^7} + \frac{5}{4}{e^5} - \frac{{43}}{{24}}{e^3}} \right){\cos ^3}E \hfill \\
   + {e^2}\left( { - \frac{1}{9}{e^6} + \frac{{17}}{{15}}{e^4} + {E^2} - \frac{{28}}{9}{e^2} + \frac{4}{3}} \right){\cos ^2}E \hfill \\
   - \frac{3}{2}e\left( { - \frac{7}{{18}}{e^6} + \frac{{37}}{{12}}{e^4} + {E^2} - \frac{{287}}{{72}}{e^2} - \frac{{119}}{{18}}} \right)\cos E \hfill \\
   - \frac{2}{9}{e^8} - \frac{{{e^6}}}{{15}} + \frac{{103}}{9}{e^4} + \left( {2{E^2} - \frac{{79}}{3}} \right){e^2} - 3{E^2} - 16 \hfill \\ 
\end{gathered}  \right)\sin E \hfill \\
   + \left( \begin{gathered}
   - {e^4}{\cos ^4}E + \left( {\frac{8}{3}{e^5} - 8{e^3}} \right){\cos ^3}E + \left( { - 8{e^4} + 29{e^2}} \right){\cos ^2}E \hfill \\
   - 16e\cos E - \frac{7}{6}{e^8} + \frac{{29}}{4}{e^6} + \frac{{27}}{8}{e^4} + \left( { - \frac{{283}}{6} + {E^2}} \right){e^2} - \frac{8}{3} - 2{E^2} \hfill \\ 
\end{gathered}  \right) \hfill \\ 
\end{gathered}  \right)\label{eqA13}
\end{equation}

\begin{equation}
{\bm{C}_{45}}\left( f \right) =0 \text{ }\text{ }\text{ }\text{ }{\bm{C}_{46}}\left( f \right) =0
\label{eqA14}
\end{equation}

\begin{equation}
 {\bm{C}_{55}}\left( f \right) =  - \frac{3}{{8{{\left( {1 - {e^2}} \right)}^{\frac{9}{2}}}}}\left( {\left( \begin{gathered}
  \frac{8}{{15}}{e^3}{\cos ^4}E - 2{e^2}{\cos ^3}E + \left( { - \frac{8}{{45}}{e^3} + \frac{8}{3}e} \right){\cos ^2}E \hfill \\
   + \left( {{e^2} - \frac{4}{3}} \right)\cos E - \frac{{16}}{{45}}{e^3} - \frac{8}{3}e \hfill \\ 
\end{gathered}  \right)\sin E + E\left( {{e^2} + \frac{4}{3}} \right)} \right)\label{eqA15}
\end{equation}

\begin{equation}
 {\bm{C}_{56}}\left( f \right) =   - \frac{{\cos E}}{{4{{\left( {1 - {e^2}} \right)}^5}}}\left( \begin{gathered}
   - \frac{4}{5}{e^3}{\cos ^4}E + \left( {{e^4} + 3{e^2}} \right){\cos ^3}E \hfill \\
   + \left( { - 4{e^3} - 4e} \right){\cos ^2}E + \left( {6{e^2} + 2} \right)\cos E - 4e \hfill \\ 
\end{gathered}  \right)\label{eqA16}
\end{equation}

\begin{equation}
{\bm{C}_{66}}\left( f \right) = \frac{9}{{4{{\left( {1 - {e^2}} \right)}^{\frac{{11}}{2}}}}}\left( \begin{gathered}
  \left( \begin{gathered}
   - \frac{4}{{45}}{e^3}{\cos ^4}E + \left( {\frac{2}{9}{e^4} + \frac{1}{3}{e^2}} \right){\cos ^3}E + \left( { - \frac{4}{{27}}{e^5} - \frac{{136}}{{135}}{e^3} - \frac{4}{9}e} \right){\cos ^2}E \hfill \\
   + \left( {{e^4} + \frac{{11}}{6}{e^2} + \frac{2}{9}} \right)\cos E - \frac{8}{{27}}{e^5} - \frac{{452}}{{135}}{e^3} - \frac{{16}}{9}e \hfill \\ 
\end{gathered}  \right)\sin E \hfill \\
   + E\left( {{e^4} + \frac{{41}}{{18}}{e^2} + \frac{2}{9}} \right) \hfill \\ 
\end{gathered}  \right)\label{eqA17}
\end{equation}

\section*{Appendix C: Flow Chart of The Numerical Method}\label{secA3}
\begin{figure}[H]
\centering
\includegraphics[width=0.78\textwidth]{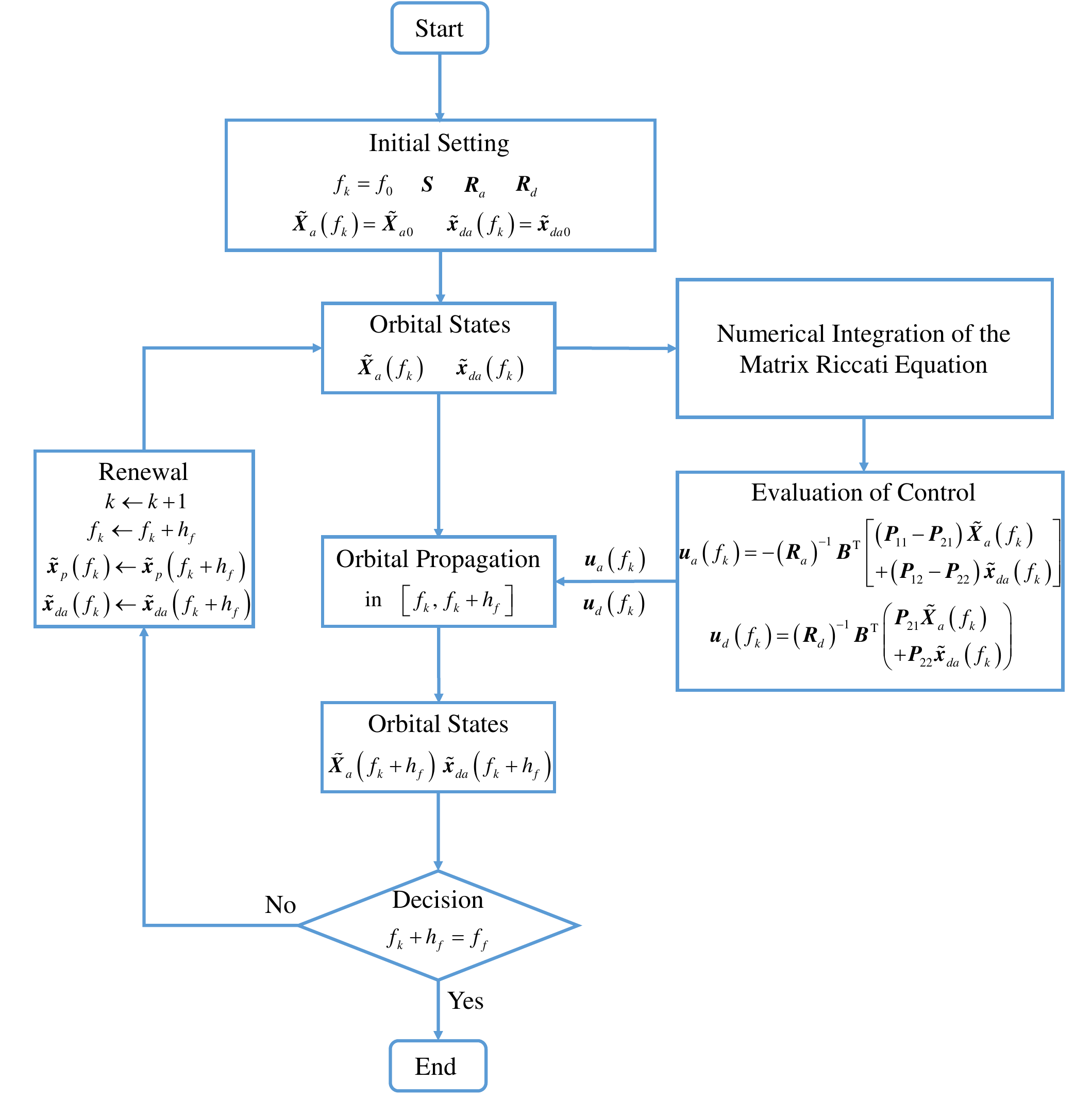}
\caption{The flow chart of the numerical method.}
\label{fig9}
\end{figure}

The flow chart of the numerical method mentioned in Section \ref{sec5} is shown in Fig. \ref{fig9}.
\section*{Funding Sources}

The fourth author acknowledges the financial support from the National Key Laboratory of Space Intelligent Control (No. HTKJ2024KL502008).

\section*{Acknowledgments}
The authors acknowledge the suggestions and assistances from Shuo Song and Yihan Peng from School of Astronautics, Beihang University.

\end{document}